\documentclass[10pt]{article}
\usepackage{amsfonts}
\usepackage{amsmath}
\usepackage{amssymb}
\usepackage{enumerate}
\usepackage{epsfig}
\usepackage{caption2}

\title{An efficient solver for problems of scattering by bodies of revolution}
\author{Youngae Han}
\newtheorem{theorem}{Theorem}[section]
\newtheorem{lemma}[theorem]{Lemma}

\newtheorem{corollary}[theorem]{Corollary}

\newtheorem{example}[theorem]{Example}
\newtheorem{definition}[theorem]{Definition}
\newenvironment{proof}{\hspace{0.5cm} {\bf Proof.}}
{$\quad {}_\blacksquare$\vspace{0.3cm}}

\begin{document}
\maketitle

\addtocounter{footnote}{1}
\begin{abstract}
This proposal relates to the design, analysis and application of a novel numerical scheme for 
the solution of axisymmetric scattering problems. To this end, a procedure is introduced to 
iteratively evaluate the solution of the Lippmann-Schwinger integral equation in $O(N\log^2 N)$ 
operations, 
where $N$ is the number of the discretization points.
The method achieves its efficiency through the use of the addition theorem and Fast Legendre Transforms (FLT).
For globally smooth sound velocities/refractive indexes the method is spectrally accurate. More generally
the order of convergence is tied to and in fact, limited by, the smoothness of the solution. 
\end{abstract}
\section{Introduction}\label{sec1}

A variety of numerical methods have been developed 
for scattering simulations. These include finite difference methods (FDM) (\cite{niraj}, \cite{teixeira}, 
\cite{guo}), 
finite element methods (FEM) (\cite{khebir}, \cite{suomin}, \cite{wolfe}), 
integral equation methods (IEM) (\cite{weiping}, \cite{zhang}, \cite{qing}) and 
methods based on fast Fourier transforms (\cite{lin}, \cite{tie}, \cite{peter}) and 
fast multipole expansions (\cite{jiming}, \cite{darve}, \cite{darrigrand}).\\

Each of these have advantages and shortcomings. For instance, FDM and FEM result 
in sparse matrices which can, therefore, be easily applied and inverted. Moreover, FDM can be implemented with 
relative ease. These methods however are typically tied to regular grids which prevents the attainment of high 
order convergence for arbitrarily shaped obstacles. In contrast, FEM, while more difficult to implement than 
its FDM counterparts can easily handle complex geometries, a feature that, in fact, largely justifies its 
popularity. A main disadvantage of the FEM approach (also shared by FDM), on the other hand, is related to the
 \emph{exterior} nature of scattering problem. Indeed, this implies that the physical \emph{radiation conditions} must be translated into conditions that allow for a finite computational domain. Although exact boundary conditions can be found \cite{mittra}, 
they are non local (in space and time). For this reason, great effort has 
gone into the design of approximate local boundary conditions (\cite{andrew}, \cite{wolfe}, \cite{jo}) 
that minimize spurious reflections. 
To achieve this, however, the conditions must typically be imposed at a certain distance from the scatterer,leading to increased computational times and memory requirements.\\

Integral equation methods, on the other hand, do not suffer from this problem, as they are based on integral 
formulas which implicitly account for radiation conditions through the use of \emph{outgoing} Green's 
functions. These methods, however, typically lead to a linear system involving  a full matrix and thus they 
are not competitive unless a specialized strategy is used to accelerate matrix-vector products. Examples of 
accelerated IEM include those based on FFTs (\cite{lin}, \cite{zhang}) 
and those that use fast multipole expansions (\cite{jiming}, \cite{kalyan}). 
Interestingly, however, to date no implementation of these exists that can attain higher than first or second
order of convergence.\\

In this paper, we present a new integral equation approach to
scattering simulations, especially designed to treat penetrable bodies
of revolution. In the spirit of methods based
on FFTs, our approach attains its efficiency through the use of fast
Legendre transforms (FLTs). In contrast with classical methods, however,
our scheme can be made to converge with arbitrarily large order. This
important additional property is achieved using ideas inspired by work in \cite{bh} and \cite{bhr} on two-dimensional scattering.
As in \cite{bh} we use the \emph{addition theorem} in the integral equation to derive a scheme whose convergence rate is tied to the global smoothness of the refractive index/sound velocity; for globally smooth 
characteristics, for instance, it converges spectrally. The contribution
of this paper is to get the high order convergence of axisymmetric case in 3-D
with two dimensional cost using FLTs.  


\section{The Helmholtz and Lippmann-Schwinger equations}\label{sec2}
As we mentioned, we are interested in the prediction of electromagnetic and acoustic scattering returns.
Helmholtz equation arises naturally in acoustics and we shall concentrate on scattering problems for 
such models. Although our results will not be directly applicable to electromagnetics, they do suggest a 
strategy to follow in this case, as the corresponding integral equations are similar in structure to those
that arise in connection with the Helmholtz model.
       
We will consider the case when the inhomogeneity is of compact support and the region under 
consideration is in $\mathbb R ^3$.
A bounded scatterer $\Omega$ is contained within a radius
$R$ from the origin. The refractive index $n(x)$ varies arbitrarily in
$\Omega$ and n(x)=1 outside the scatterer. Thus, setting
m(x)=$1-n(x)^2$ we have m(x)=0 for x outside $\Omega$
and a scattered field $u^s$ will be generated by an incident field $u^i$ satisfying 
the linearized equations in free space

\begin{equation}\label{1}
 \Delta {u^{i}}+k^2{u^{i}} = 0 \mbox{ in } {\mathbb R ^3}.
\end{equation}

The total field $ {u}={u^{i}}+{u^{s}}$ will then satisfy

\begin{equation}\label{2}
 \Delta {u}+k^2n^2{u} = 0 \mbox{ in } {\mathbb R ^3}
\end{equation}
with appropriate conditions at infinity. The physically relevant condition is
the \emph{Sommerfeld radiation condition} which guarantees that
the scattered field is outgoing,
\begin{equation}\label{3}
 \lim\limits_{r \rightarrow \infty} r\, \bigg (\frac{\partial u^{s} }{\partial r}-ik u^{s} \bigg)=0. 
\end{equation}

In fact, as we mentioned, this last condition constitutes one of the central challenges in the design of 
accurate numerical methods, especially if it must be approximated (see e.g., \cite{teixeira}, \cite{kyung}).
Such approximations, however, can be avoided in an integral formulation of the problem that 
explicitly accounts for the outgoing character of the solution. 
Indeed, suppose $u\in C^2( \mathbb R ^3 )$ is a solution of (\ref{2})-(\ref{3}) and let 
\begin{equation}\label{4}
\Phi(x,y)=\frac{1}{4\pi}\frac{e^{ik\mid x-y\mid}}{\mid x-y\mid},
\end{equation}
denote the free-space outgoing Green's function. Let $x \in \mathbb R ^3 $
be an  arbitrary point and choose an open ball $B$ of radius r with exterior  unit normal $\nu$ 
containing the support of $m$ and such that $x \in B$. Then, we arrive at the 
$\textit{Lippmann-Schwinger equation}$ which is equivalent to (\ref{2}) and (\ref{3}), 
\begin{equation}\label{5}
u(x)= u^{i}(x)-k^{2}\int\limits_{B}^{}\! \Phi(x,y) u(y)m (y)\,dy,\ x \in \mathbb R^3
\end{equation}
where $B$ is any ball containing the support of $m$.

\section{Axisymmetric formulation of the Lippmann-Schwinger integral equation}\label{sec3}
A convenient form of the Lippmann-Schwinger equation can be
derived if we appeal to the addition theorem.
The addition theorem (see \cite{colton}) states that 
\begin{equation}\label{6}
\begin{split}
\frac{1}{4\pi}\frac{e^{ik\mid x-y\mid}}{\mid x-y\mid}
&=ik\sum_{n=0}^{\infty}h_n^{(1)}(k\mid \widetilde{\rho}\mid) j_{n}(k\mid \rho \mid) 
\frac{2n+1}{4\pi}\bigg[P_n(\cos\widetilde{\theta})P_n(\cos\theta)\\
&+2\sum_{m=1}^n
\frac{(n-m)!}{(n+m)!}P_n^m(\cos\widetilde{\theta})P_n^m(\cos\theta)\cos(m(\widetilde{\phi}-\phi))\bigg]
\end{split}
\end{equation}
where $x=(\widetilde{\rho}\sin \widetilde{\theta}\cos \widetilde{\varphi},
\widetilde{\rho} \sin \widetilde{\theta} \sin \widetilde{\varphi},
\widetilde{\rho} \cos \widetilde{\theta} ),$ $y=(\rho \sin \theta \cos
\varphi,\rho \sin \theta \sin \varphi,\rho \cos \theta ),$ and  $P_{n}^m$ are the
associated Legendre functions ($P_{n}^0=P_{n}$ are the Legendre 
polynomials). 
The series and its term by term first derivatives with
respect to $\mid x\mid$ and $\mid y\mid$ are absolutely convergent on
compact subsets of $\mid x \mid > \mid y \mid$.\\
Therefore, we obtain, using the addition theorem
and lemma, that 
\begin{equation}\label{7}
\begin{split}
  \int\limits_{\mathbb R^3}^{}\! \Phi(x,y) u(y)m (y)\,dy
  &= \int\limits_{\mathbb R^3}^{}\! 
        \frac{1}{4\pi}\frac{e^{ik\mid x-y\mid}}{\mid x-y\mid} u(y)m (y)\,dy\\
  & =ik\int\limits_{\mathbb R^3}^{}\! \sum_{n=0}^{\infty}
  h_{n}^{(1)}(k\rho_{>})
  j_{n}(k\rho_{<})\frac{2n+1}{4\pi}\bigg[P_n(\cos\widetilde{\theta})P_n(\cos\theta)\\
 &+2\sum_{m=1}^n\frac{(n-m)!}{(n+m)!}P_n^m(\cos\widetilde{\theta})P_n^m(\cos\theta)\cos(m(\widetilde{\phi}-\phi))\bigg]\\
 &   \cdot u\cdot m \cdot \rho^2 \sin\theta\,d\rho d\theta d\varphi,\\
\end{split}
\end{equation}
where
$$
 \rho_{>}=\mbox{max}(\widetilde{\rho},\rho),\\ \rho_{<}=\min(\widetilde{\rho},\rho).
$$

If the incoming ray and the obstacle are axisymmetric then, by
uniqueness \cite{colton}, the solution $u$ of (\ref{2}) and (\ref{3})
will be axisymmetric.
Then (\ref{7}) simplifies to
\begin{equation}\label{8}
\begin{split}
 \int\limits_{\mathbb R^3}^{}\! \Phi(x,y) u(y)m (y)\,dy 
 &=ik\sum_{n=0}^{\infty} \frac{(2n+1)}{2}
\int\limits_{0}^{\pi}\int\limits_{0}^R\!  
  h_{n}^{(1)}(k\rho_{>}) j_{n}(k\rho_{<})P_{n}(\cos\theta) P_{n}(\cos \widetilde{\theta})\\
&\cdot u(\rho,\theta) m (\rho,\theta)\rho^2\sin\theta \,d\rho d\theta.
\end{split}
\end{equation}

Finally, expanding $u$ and $u^i$ in Legendre series we obtain

\begin{multline}\label{9}
\sum_{n=0}^{\infty}u_n(\widetilde{\rho}) P_{n}(\cos \widetilde{\theta})
=\sum_{n=0}^{\infty}u_n^i(\widetilde{\rho}) P_{n}(\cos \widetilde{\theta})\\
   -\sum_{n=0}^{\infty}\frac{(2n+1)ik^3}{2} \int\limits_{0}^{\pi}\int\limits_{0}^R\!
  h_{n}^{(1)}(k\rho_{>}) j_{n}(k\rho_{<})   P_{n}(\cos\theta)\\ 
\cdot u(\rho,\theta) m (\rho,\theta)\rho^2\sin\theta 
   \,d\rho d\theta P_{n}(\cos \widetilde{\theta}),
\end{multline}
from which the axisymmetric formulation 
\begin{equation}\label{10}
u_n(\widetilde{\rho})=u_n^i(\widetilde{\rho})
-\frac{(2n+1)ik^3}{2} \int\limits_{0}^{\pi}\int\limits_{0}^R\!  
h_{n}^{(1)}(k\rho_{>}) j_{n}(k\rho_{<})
 P_{n}(\cos\theta)u(\rho,\theta) m (\rho,\theta)\rho^2\sin\theta 
\,d\rho d\theta
\end{equation}
follows.

\section{Numerical implementation}
To implement (\ref{10}) numerically,  we approximate 
$u$ and $m$ by a truncated Legendre series
$$
  u^{F}(\rho,\cos\theta) = \sum_{n=0}^{F} u_{n}(\rho) P_{n}(\cos\theta),
$$
$$
  m^{M}(\rho,\cos\theta) = \sum_{n=0}^{M} u_{n}(\rho) P_{n}(\cos\theta).
$$

Then, if $n \le F$, using the orthogonality properties of the Legendre polynomials we have
\begin{multline}\label{11}
  \int\limits_{0}^{\pi}\!u^F(\rho,\cos \theta)P_n(\cos \theta)m(\rho,\cos \theta)\,\sin\theta d\theta\\
=\int\limits_{0}^{\pi}\!u^F(\rho,\cos \theta)P_n(\cos \theta)\sum_{l=0}^{2F} 
m _{l}(\rho) P_{l}(\cos \theta)\,\sin\theta d\theta ,
\end{multline}

so that we $m$ can be approximated by 
$$m ^{2F}(\rho,\cos\theta) = \sum_{l=0}^{2F} m _{l}(\rho) P_{l}(\cos\theta).$$

The formulation (\ref{10}) thus reduces to 
\begin{equation}\label{12}
u_n(\widetilde{\rho})=u_n^i(\widetilde{\rho})
+\frac{i}{2}K_{n}(\widetilde{\rho}),  
\end{equation}
where
\begin{equation}\label{13}
 K_{n}(\widetilde{\rho})=-(2n+1)k^3\int\limits_{0}^{\mathbf R}\!
h_{n}^{(1)}(k\rho_{>}) j_{n}(k\rho_{<}) I_{n}(\rho)\rho^2\,d\rho,
\end{equation}
and
\begin{equation}\label{14}
 I_{n}(\rho)
=\int\limits_{0}^{\pi}\!u^{F}(\rho,\theta) m ^{2F}(\rho,\theta) 
P_{n}(\cos\theta)\sin\theta\,d\theta,
\end{equation}
$n=0,1,\dots,F$. To solve (\ref{12}), we shall appeal to an iterative method, e.g., GMRES, which 
reduces the overall problem to the fast evaluation of $K_{n}(\widetilde{\rho})$ in (\ref{12}).
To achieve the latter we propose below an efficient strategy to calculate the angular 
integration in (\ref{14}) and the subsequent radial integration in (\ref{13}).

\subsection{Angular integration}
A natural approach to the evaluation of (\ref{14}) is to exploit the orthogonality of the Legendre
 polynomials by expanding the function $u^Fm^{2F}(\rho,\theta)$ in Legendre series for fixed $\rho$. This, 
of course, demands that we calculate the \emph{Legendre transform} $\{c_n\}$
\begin{equation}\label{15}
   u^F(\rho,\theta)m^{2F}(\rho,\theta)=\sum_{n=0}^{3F}c_n(\rho)P_n(\cos\theta),
\end{equation}
which, in principle, requires $O(F^2)$ operations (for fixed $\rho$).
Fortunately, however, algorithms have been designed for the fast evaluation of the map   
\begin{equation}\label{16}
   \{f(\theta_j)\}_{j=1}^F \rightarrow \{c_n\}_{n=1}^F,
\end{equation}
where
$$
 f(\theta_j)=\sum_{n=0}^{F}c_nP_n(\cos\theta_j),
$$
and its inverse, for appropriate choices of the angles $\{\theta_j\}$; see Appendix.
These algorithms, which we shall refer to as \emph{Fast Legendre Transform} (FLT) and  
\emph{Inverse Fast Legendre Transform} (IFLT), then suggest an implementation of the angular integration 
in the form

\begin{equation}\label{17}
\{  I_{n}(\rho)/\tau_{n}\}_{n=0}^{F}=FLT_{3F}(IFLT_{3F}(\{u_{n}(\rho)\}_{n=0}^{F})\cdot
IFLT_{3F}(\{m_{n}(\rho)\}_{n=0}^{2F})),
\end{equation}
where
$$\tau_{n}=\int\limits_{-1}^{1}\!{\mid p_{n}(x)\mid}^2\,dx=\frac{2}{2n+1}.$$

\subsection{Radial integration}
As can be easily checked, the integrand in (\ref{13}) has a corner-type singularity at
$\rho=\widetilde{\rho},$ which suggests that we write 

\begin{equation}\label{18}
\begin{split}
 \frac{-K_{n}(a)}{(2n+1)k^3}
 & =h_{n}^{(1)}(ka)\int\limits_{0}^{\min(a,\mathbf R)}\!j_{n}(k\rho)
   I_{n}(\rho)\rho^2\,d\rho
  +j_{n}(ka)\int\limits_{\min(a,\mathbf R)}^{\mathbf R}\!h_{n}^{(1)}(k\rho) 
   I_{n}(\rho)\rho^2\,d\rho \\
 &=i\bigg[ y_{n}(ka)\int\limits_{0}^{\min(a,\mathbf R)}\!j_{n}(k\rho)
   I_{n}(\rho)\rho^2\,d\rho
  +j_{n}(ka)\int\limits_{\min(a,\mathbf R)}^{\mathbf R}\!y_{n}(k\rho) 
   I_{n}(\rho)\rho^2\,d\rho \bigg] \\
 & +j_{n}(ka)\int\limits_{0}^{\mathbf R}\!j_{n}(k\rho)I_{n}(\rho)\rho^2\,d\rho. \\
\end{split} 
\end{equation}
Since, $h_n^{(1)}(k\rho)$ blows up at the origin although
$j_{n}(k\rho)h_{n}^{(1)}(k\rho)$ behaves nicely. 
So, we define the modified Bessel functions 
$\widetilde{j_n}(\rho),\mbox{ }\widetilde{y_n}(\rho)$ as follows.
 
\begin{equation}\label{19}
\begin{split}
&\widetilde{j_n}(\rho):=\frac{1\cdot3\cdot5\cdot\dots(2n+1)}{\rho^n}j_n(\rho)=
\bigg[1-\frac{\frac{1}{2}\rho^2}{1!(2n+3)}+\frac{({\frac{1}{2}\rho^2})^2}{2!(2n+3)(2n+5)}+\dots
\bigg]\\
&\widetilde{y_n}(\rho):=\frac{\rho^{n+1}}{-1\cdot1\cdot3\cdot5\cdot\dots(2n-1)}y_n(\rho)=
\bigg[1-\frac{\frac{1}{2}\rho^2}{1!(1-2n)}+\frac{({\frac{1}{2}\rho^2})^2}{2!(1-2n)(3-2n)}+\dots
\bigg]\\
\end{split}
\end{equation}

Then, the equation (\ref{18}) becomes
\begin{equation}\label{20}
\begin{split}
K_{n}(a)
 &=i\bigg[ \widetilde{y_n}(ka)\int\limits_{0}^{\min(a,\mathbf R)}\!
\big(\frac{\rho}{a}\big)^{n+1}\widetilde{j_n}(k\rho)
   I_{n}(\rho)k^2\rho\,d\rho \\
 &+\widetilde{j_n}(ka)\int\limits_{\min(a,\mathbf R)}^{\mathbf R}\!
\big(\frac{a}{\rho}\big)^{n}\widetilde{y_n}(k\rho) 
   I_{n}(\rho)k^2\rho\,d\rho \bigg] \\
 & +\widetilde{j_n}(ka)(ka)^n\int\limits_{0}^{\mathbf R}\!
\frac{(k\rho)^n(-2n-1)k^3}{1\cdot3^2\cdot5^2\dots (2n+1)^2}\widetilde{j_n}(k\rho)I_{n}(\rho)\rho^2\,d\rho. \\
\end{split} 
\end{equation}
and approximate each integrand in the right-hand side separately. To this end, 
we divide the integration domain in a number $N_{i}$ of equi-length
interpolation intervals $U_{j}=[u_j^0,u_j^1]$, $1\le j \le N_i$ on which we approximate 

\begin{equation}\label{21}
I_{n}(\rho)\thickapprox\sum_{m=0}^{N_d-1}c_{m}^j T_{m}^{u_j^0,\mbox{}u_j^1}(\rho),
\end{equation}
for $\rho \in U_j,$ where
$$
 T_{m}^{u_j^0,\mbox{}u_j^1}(\rho)=  T_{m}\bigg(\frac{\rho-(u_j^1+u_j^0)/2}{(u_j^1-u_j^0)/2}\bigg)
$$
are Chebyshev polynomials in $U_{j}$.
Then letting $\{\rho^j_k\}^{N_d}_{k=1}$ denote the Chebyshev points in $U_{j}$ we have 
\begin{equation}\label{22}
\begin{split}
 K_{n}(\rho^j_{k}) &
 \thickapprox   i\widetilde{y_{n}}(k\rho^j_{k})  \bigg[\sum_{p=1}^{j-1}
    \sum_{m=0}^{N_d-1}c_{m}^{p}\int\limits_{U_p}\!\big(\frac{\rho}{\rho^j_k}\big)^{n+1}
\widetilde{j_{n}}(k\rho)
      T_{m}^{u_p^0,\mbox{}u_p^1}(\rho)k^2\rho\,d\rho    \\
 & +\sum_{m=0}^{N_d-1}c_{m}^{j}\int\limits_{u_j^0}^{\rho^j_{k}}
 \!\big(\frac{\rho}{\rho^j_k}\big)^{n+1}\widetilde{j_{n}}(k\rho)
   T_{m}^{u_j^0,\mbox{}u_j^1}(\rho)k^2\rho\,d\rho \bigg]   \\
 & +i\widetilde{j_{n}}(k\rho^j_{k})  \bigg[ \sum_{p=j+1}^{N_i}
\sum_{m=0}^{N_d-1}c_{m}^{p}\int\limits_{U_p}\!\big(\frac{\rho^j_k}{\rho}\big)^{n}
\widetilde{y_{n}}(k\rho)
 T_{m}^{u_p^0,\mbox{}u_p^1}(\rho)k^2\rho\,d\rho  \\
 & +\sum_{m=0}^{N_d-1}c_{m}^{j}\int\limits^{u_j^1}_{\rho^j_{k}}
 \!\big(\frac{\rho^j_k}{\rho}\big)^{n}
\widetilde{y_{n}}(k\rho)
 T_{m}^{u_j^0,\mbox{}u_j^1}(\rho)k^2\rho\,d\rho \bigg]\\
 &
 +\widetilde{j_n}(k\rho^j_{k})(k\rho^j_{k})^n\sum_{p=1}^{N_i}\sum_{m=0}^{N_d-1}c_{m}^{p}
\int\limits_{U_p}\!\frac{(k\rho)^n(-2n-1)k^3}{1\cdot3^2\cdot5^2\dots
  (2n+1)^2}\widetilde{j_n}(k\rho)T_{m}^{u_p^0,\mbox{}u_p^1}(\rho)\rho^2\,d\rho, \\
\end{split}
\end{equation}

which can be readily evaluated if the \emph{moments} 
\begin{equation*}
\begin{split}
&\int\limits_{\rho^j_{k}}^{\rho^j_{k+1}}\!\big(\frac{\rho}{\rho^j_{k+1}}\big)^{n+1}
\widetilde{j_{n}}(k\rho)
 T_{m}^{u_j^0,\mbox{}u_j^1}(\rho)k^2\rho\,d\rho
\mbox{ ,  }
\int\limits_{\rho^j_{k}}^{\rho^j_{k+1}}\!\big(\frac{\rho^j_k}{\rho}\big)^{n}
\widetilde{y_{n}}(k\rho)
 T_{m}^{u_j^0,\mbox{}u_j^1}(\rho)k^2\rho\,d\rho,\\
&\int\limits_{\rho^j_{k}}^{\rho^j_{k+1}}\!\frac{(k\rho)^n}{1\cdot3\cdot5\cdots(2n+1)}
\widetilde{j_{n}}(k\rho)
 T_{m}^{u_j^0,\mbox{}u_j^1}(\rho)\rho^2\,d\rho
\end{split}
\end{equation*}

$(\rho^j_{0}=u_j^0 \mbox{, }\rho^j_{(N_d+1)}=u_j^1)$
are pre-computed and stored for 
$$
0 \leq n \leq F\mbox{, }0 \leq m \leq N_d-1\mbox{, } 1 \leq j \leq N_i \mbox{, } 0 \leq k \leq N_d.
$$

\subsection{Algorithm and operation count}
The prescriptions in \S 4.1 and \S 4.2 can be summarized in the following algorithm.\\
\emph{
\textbf{Algorithm}\\
\textbf{Input} \mbox{  }
$\{u_n^i(\rho_k^j)\}_{n=0}^{F}$\mbox{ ,  }$\{m_n(\rho_k^j)\}_{n=0}^{2F}$ : 
Legendre transform coefficients of $u^i(\rho_k^j,\theta)$ and $m(\rho_k^j,\theta).$\\
\textbf{Output} 
$\{u_n(\rho_k^j)\}_{n=0}^{F}$ : 
Legendre transform coefficients of $u(\rho_k^j,\theta).$\\
\textbf{Stages}\\
0 stage,\\
\begin{itemize}
\item compute Bessel functions
$$
 \widetilde{y_{n}}(\emph{k}\rho^j_{k})\mbox{ ,  } \widetilde{j_{n}}(\emph{k}\rho^j_{k})
$$
\item compute ratios
$$
  \frac{(\rho^j_k)^{n+1}}{(\rho^j_{k+1})^{n+1}}\mbox{ ,  }
\frac{(\rho^j_{k})^{n}}{(\rho^j_{k+1})^{n}}
$$
\item compute moments $\alpha_{j,k,n,m}$, $\beta_{j,k,n,m}$, $\gamma_{j,k,n,m}$,    
\begin{equation*}
\begin{split}
&\alpha_{j,k,n,m}:=\int\limits_{\rho^j_{k}}^{\rho^j_{k+1}}\!\big(\frac{\rho}{\rho^j_{k+1}}\big)
^{n+1}\widetilde{j_{n}}(k\rho)T_{m}^{u_j^0,\mbox{}u_j^1}(\rho)k^2\rho\,d\rho
\mbox{ ,  }
\beta_{j,k,n,m}:=\int\limits_{\rho^j_{k}}^{\rho^j_{k+1}}\!\big(\frac{\rho^j_k}{\rho}\big)^{n}
\widetilde{y_{n}}(k\rho)T_{m}^{u_j^0,\mbox{}u_j^1}(\rho)k^2\rho\,d\rho,\\
&\gamma_{j,k,n,m}:=\int\limits_{\rho^j_{k}}^{\rho^j_{k+1}}\!\frac{(k\rho)^n}{1\cdot3\cdot5\cdots(2n+1)}
\widetilde{j_{n}}(k\rho)
 T_{m}^{u_j^0,\mbox{}u_j^1}(\rho)\rho^2\,d\rho\\
\end{split}
\end{equation*}
(e.g. with a $Clenshaw-Curtis\mbox{  } quadrature$ \cite{recipe})
\item set  $u_n(\rho_k^j)=u_n^i(\rho_k^j)$\\
\end{itemize}
k stage,$\mbox{ for k=1 to the number of iterations required for GMRES}$
$\mbox{ to converge to a prescribed }$\\${tolerance}\ $\\
\begin{itemize}
\item compute angular integration $\{I_n(\rho^j_{k})\}_{n=0}^{F}$ 
for given $u_n$,
$$
\{  I_{n}(\rho^j_{k})/\tau_{n}\}_{n=0}^{F}=FLT_{3F}(IFLT_{3F}(\{u_{n}(\rho^j_{k})\}_{n=0}^{F})\cdot
IFLT_{3F}(\{m_{n}(\rho^j_{k})\}_{n=0}^{2F}))
$$
\item compute radial integration
\renewcommand{\labelenumi}{{\rm(\alph{enumi})}}
\begin{enumerate}
\item  compute $\{c_{m}^j\}_{m=0}^{N_d-1}$
$$
I_{n}(\rho)\thickapprox\sum_{m=0}^{N_d-1}c_{m}^j T_{m}^{u_j^0,\mbox{}u_j^1}(\rho) \mbox{ for }
\rho \in U_j \ ;
$$
\item define $\mu_{j,k,n}$, $\zeta_{j,k,n}$,  $\xi_{j,k,n}$
$$
\mu_{j,k,n}:=\sum_{m=0}^{N_d-1}c_{m}^j\alpha_{j,k,n,m}
\mbox{ ,  }
\zeta_{j,k,n}:=\sum_{m=0}^{N_d-1}c_{m}^j\beta_{j,k,n,m}
\mbox{ ,  }
\xi_{j,k,n}:=\sum_{m=0}^{N_d-1}c_{m}^j\gamma_{j,k,n,m} ;
$$
\item compute $K_{n}(\rho^j_{k})$
\begin{equation*}
\begin{split}
&\aleph=\sum_{p=1}^{N_i}\sum_{l=0}^{N_d}\xi_{p,l,n}\\ 
&s(\rho_1^1)=\widetilde{y_n}(k\rho_1^1)\mu_{1,0,n}\\
&s(\rho_2^1)=\frac{\widetilde{y_n}(k\rho_2^1)}{\widetilde{y_n}(k\rho_1^1)}\frac{(\rho^1_1)^{n+1}}
{(\rho_2^1)^{n+1}}s(\rho_1^1)+\widetilde{y_n}(k\rho_2^1)\mu_{1,1,n}\\
&s(\rho_3^1)=\frac{\widetilde{y_n}(k\rho_3^1)}{\widetilde{y_n}(k\rho_2^1)}\frac{(\rho^1_2)^{n+1}}
{(\rho_3^1)^{n+1}}s(\rho_2^1)+\widetilde{y_n}(k\rho_3^1)\mu_{1,2,n}\\
&\dots\\
&s(\rho_{N_d}^{N_i})=\frac{\widetilde{y_n}(k\rho_{N_d}^{N_i})}{\widetilde{y_n}(k\rho_{N_d-1}^{N_i})}
\frac{(\rho_{N_d-1}^{N_i})^{n+1}}{(\rho_{N_d}^{N_i})^{n+1}}s(\rho_{N_d-1}^{N_i})+
\widetilde{y_n}(k\rho_{N_d}^{N_i})\mu_{N_i,N_d-1,n}\\
& \\
&q(\rho_{N_d}^{N_i})=\widetilde{j_n}(k\rho_{N_d}^{N_i})\zeta_{N_i,N_d,n}\\
&q(\rho_{N_d-1}^{N_i})=\frac{\widetilde{j_n}(k\rho_{N_d-1}^{N_i})}
{\widetilde{j_n}(k\rho_{N_d}^{N_i})}\frac{(\rho_{N_d-1}^{N_i})^{n}}{(\rho_{N_d}^{N_i})^{n}}
q(\rho_{N_d}^{N_i})+\widetilde{j_n}(k\rho_{N_d-1}^{N_i})\zeta_{N_i,N_d-1,n}\\
&q(\rho_{N_d-2}^{N_i})=\frac{\widetilde{j_n}(k\rho_{N_d-2}^{N_i})}
{\widetilde{j_n}(k\rho_{N_d-1}^{N_i})}\frac{(\rho_{N_d-2}^{N_i})^{n}}{(\rho_{N_d-1}^{N_i})^{n}}
q(\rho_{N_d-1}^{N_i})+\widetilde{j_n}(k\rho_{N_d-2}^{N_i})\zeta_{N_i,N_d-2,n}\\
&\cdots\\
&q(\rho_1^1)=\frac{\widetilde{j_n}(k\rho_{1}^{1})}
{\widetilde{j_n}(k\rho_{2}^{1})}\frac{(\rho_{1}^{1})^{n}}{(\rho_{2}^{1})^{n}}
q(\rho_{2}^{1})+\widetilde{j_n}(k\rho_{1}^{1})\zeta_{1,1,n}\\
&\\
&K_{n}(\rho^j_{k}):=i(s(\rho^j_{k})+q(\rho^j_{k}))-k^3\frac{(k\rho^j_{k})^n}{1\cdot3\cdot5
\cdots(2n-1)}\widetilde{j_n}(k\rho^j_{k})\aleph
\end{split}
\end{equation*}
\item define
$$
u_n(\rho^j_k):=u_n(\rho^j_k)
-\frac{i}{2}K_{n}(\rho^j_k).
$$ 
\end{enumerate}
\end{itemize}
}
As follows from the description above, at each iteration, the number of multiplications is given by
\begin{center}
\renewcommand{\arraystretch}{1.25}
\begin{tabular}{|p{1.6in}| p{1.6in}|}
\hline
\textbf{Stage}  & \textbf{number of operations}\\ \hline
$\mbox{ Angular integration }$ & $O(F(\log F)^2N_dN_i)$ \\ \hline 
$\mbox{ Radial integration } $ & $(a)+(b)+(c)$ \\ \hline
 \mbox{  }  (a)  & $O(N_d (\log N_d)N_iF)$ \\ \hline
 \mbox{  }  (b)  & $O(N_d^2N_iF)$ \\ \hline
 \mbox{  }  (c)  & $O(N_iN_dF)$ \\ \hline
\end{tabular}
\end{center}

Therefore, the total number of operations per iteration is
$$O\big[N_iN_d F ( \log^2 F+\log N_d+ N_d+1)\big] = O(N \log^2 N)$$
where $N=N_iN_d F$ is the total number of degrees of freedom.

\section{Numerical results}
In this section we present results that confirm the predicted performance of the algorithms described in $\S$ 4.3.

\begin{example}
In our first example we consider the scattering off a homogeneous sphere of radius 1 and index $n(x)=2$.
In this case the problem is explicitly solvable, so that comparison with the exact solution is possible.
Indeed, for a plane wave incident in the positive $z$-direction we have

\begin{equation}\label{23}
u^{i}=e^{ik\vec{x}(0,0,1)}=e^{ik\rho \cos\theta}=\sum_{n=0}^{\infty}i^n(2n+1)j_{n}(k\rho)
 P_n(\cos\theta).
\end{equation}

Equation (\ref{23}), together with the Helmholtz equation  (\ref{2}), 
and the radiation condition (\ref{3}) readily deliver
$$
  u^{s}=
\begin{cases}
\sum_{n=0}^{\infty}\{a_{n}j_{n}(2k\rho)-i^n(2n+1)j_{n}(k\rho)\} P_n(\cos\theta), 
&\rho \le1\\
\sum_{n=0}^{\infty}b_{n}h_{n}(k\rho) P_n(\cos\theta), &\rho \ge 1\\
\end{cases}
$$

where $a_n$ and $b_n$ solve
\begin{equation*}
  \left(
\begin{matrix}
-j_{n}(2k) & h_{n}(k) \\
-2j_{n}'(2k) & h_{n}'(k)
\end{matrix}
\right)
\left(
\begin{matrix}
a_n \\
b_n
\end{matrix}
\right)=-
\left(
\begin{matrix}
j_{n}(k)\\
j_{n}'(k)
\end{matrix}
\right)i^n(2n+1).
\end{equation*}
Therefore,

$$
  u^{s}=
\begin{cases}
\sum_{n=0}^{\infty}\{ \widetilde{a_{n}}\rho^n\widetilde{j_{n}}(2k\rho)
-i^n\frac{(k\rho)^n}{1\cdot3\cdot\cdots(2n-1)}\widetilde{j_{n}}(k\rho)\} P_n(\cos\theta), 
&\rho \le1\\
\sum_{n=0}^{\infty}\big(\frac{\widetilde{b}_{n}k^{2n+1}}{-1\cdot1^2\cdot3^2\cdots(2n-1)^2(2n+1)}
\rho^n\widetilde{j_{n}}(k\rho)+i\widetilde{b_n}\frac{1}{\rho^{n+1}}\widetilde{y_n}(k\rho)\big)P_n(\cos\theta), &\rho \ge 1\\
\end{cases}
$$

where $\widetilde{a_{n}}$ and $\widetilde{b_n}$ solve
\begin{equation*}
\begin{split}
  \left(
\begin{matrix}
\widetilde{j_{n}}(2k),& -\delta_n\widetilde{j_{n}}(k)-i\widetilde{y_n}(k) \\
n\widetilde{j_{n}}(2k)+2k\widetilde{j_{n}}'(2k),&
-\delta_n(n\widetilde{j_{n}}(k)+\widetilde{j_{n}}'(k)k)+i(n+1)\widetilde{y_n}(k)-i\widetilde{y_n}'(k)k
\end{matrix}
\right)&
\left(
\begin{matrix}
\widetilde{a_{n}}\\
\widetilde{b_n}
\end{matrix}
\right)\\&=
\left(
\begin{matrix}
\frac{i^nk^n}{\varpi_n}\widetilde{j_{n}}(k)\\
\frac{i^nk^n}{\varpi_n}(n\widetilde{j_{n}}(k)+\widetilde{j_{n}}'(k)k)
\end{matrix}
\right).
\end{split}
\end{equation*}
where
$$
  \delta_n:=\frac{k^{2n+1}}{-1\cdot1^2\cdot3^2\cdots(2n-1)^2(2n+1)}
$$
$$
  \varpi_n:=
\begin{cases}
1\cdot3\cdot\cdots(2n-1), 
&n \ge1\\
1&n=0
\end{cases}
$$
In figure 1 we show the error $E=\|u_{exact}-u_{approx}\|_{\infty}$ between the approximate and exact 
solutions for different values of 
the interpolation orders $N_d$ in (\ref{21}).The figure demonstrates that, indeed, the radial integrator 
incurs an error of  $O((\frac{1}{N_iN_d})^{N_d})$. Numerical values corresponding to this figure are 
displayed in table 1 where we also exhibit timings and iteration counts.
\end{example}

\epsfxsize=7 cm
\epsfysize=5 cm
\begin{center}
\epsfbox{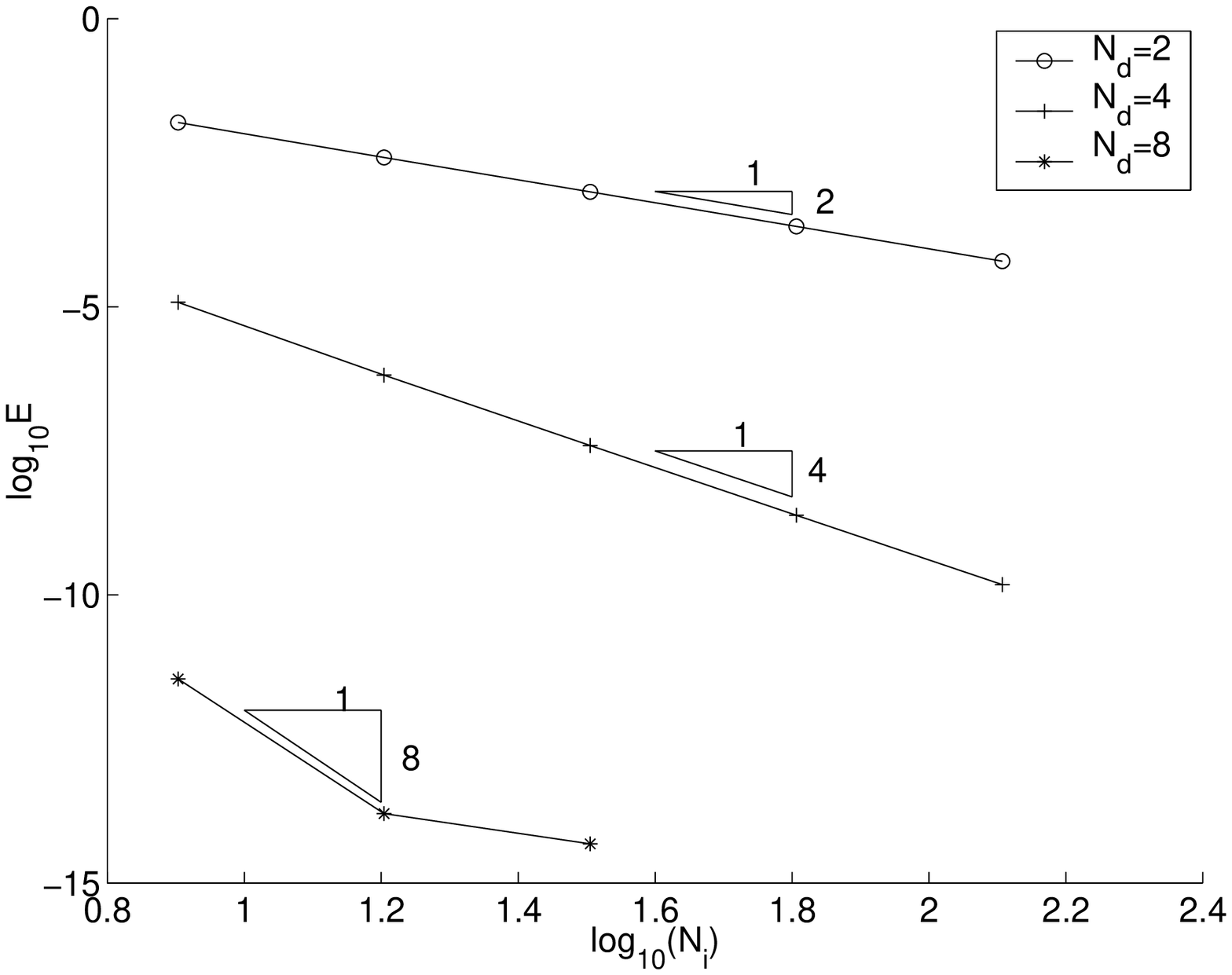}
\end{center}
$$\mbox{ Figure 1: Error in radial integrator for Example 5.1 (sphere). }$$\\

\begin{center}
\renewcommand{\arraystretch}{1.25}
\begin{tabular}{|r|p{1.6in} |r  |r|}
\hline
\textbf{$N_i$}  & \textbf{time per iteration}&\textbf{error}&\textbf{$\log_2$(error ratio)}\\ \hline
$2^3$ & 3  (sec)  & 1.19372e-05 & \\ \hline
$2^4$ & 8  (sec) & 6.51158e-07 & 4.1963  \\ \hline
$2^5$ & 15  (sec) & 3.91409e-08 & 4.0563 \\ \hline
$2^6$ & 30  (sec) & 2.41378e-09 & 4.0193\\ \hline
$2^7$ & 61 (sec)&  1.50428e-10 & 4.0042\\ \hline
\end{tabular}
\end{center}
 $$\mbox{ Table 1: Results corresponding to Example 5.1: sphere.}$$ 
 $$\mbox{ Parameters: } F=2^8-1,\ N_d=4,\ 0\le \rho \le 4.$$

\begin{example}
To test the angular integration, our second example relates the sphere
centered at $(0,0,d)$ with radius $r$. $u_i(\rho,\theta)$ is given as follows,
\begin{equation}\label{24}
\begin{split}
u^{i}&=e^{ik(\rho\cos\theta-d)}\\
&=e^{-ikd}\sum_{n=0}^{\infty}i^n(2n+1)j_{n}(k\rho)
 P_n(\cos\theta).\\
\end{split}
\end{equation}


Since,
\begin{equation*}
\begin{split}
\int\limits_{0}^{\theta_0}\!P_m(\cos\theta)\sin\theta\,d\theta&=
\int\limits_{\cos\theta_0}^1\!P_m(t)\,dt\\
&=\sqrt{(1-\cos^2\theta_0)}P_m^{-1}(\cos\theta_0)\\
&=\begin{cases}
\sqrt{(1-\cos^2\theta_0)}\frac{(m-1)!}{(m+1)!}P_m^1(\cos\theta_0), 
&m \ge1\\
(1-\cos\theta_0), &m=0.
\end{cases}\\
\end{split}
\end{equation*}
where
$$
  P_n^m=(1-x^2)^{\frac{m}{2}}\frac{d^mP_n(x)}{dx^m}.
$$
The exact Legendre series coefficients of $m(\rho,\theta)$ is given as follows.
if  $0 \le(d-r) \le\rho \le(d+r)$, 
\begin{equation*}
\begin{split}
m(\rho,&\theta)=\\
&(1-n_0^2)\Big(\sum_{m=1}^{\infty}\frac{2m+1}{2}\sqrt{(1-\cos^2\theta_0)}\frac{(m-1)!}{(m+1)!}
P_m^1(\cos\theta_0)P_m(\cos\theta)+ \frac{1}{2}(1-\cos\theta_0)P_0(\cos\theta)\Big)\\
\end{split}
\end{equation*}
where $(\rho\sin\theta_0,\mbox{  }\rho\cos\theta_0)$ is a solution of 
$$
y^2+z^2=\rho^2,\mbox{  }y^2+(z-d)^2=r^2.
$$ 
So, we have the exact Legendre coefficients of $u$ and $m$ and an approximated
solution is compared with the exact solution from example 1. Tables 2-5
shows the radial and angular convergences of a discontinuous scatterer. 
\end{example}

\begin{center}
\renewcommand{\arraystretch}{1.25}
\begin{tabular}{|r|p{1.6in} |r  |r|}
\hline
\textbf{$N_i$} &\textbf{time per iteration}  &\textbf{error}&\textbf{error ratio}\\ \hline
$2^3$&1  (sec) &0.0650681  & \\ \hline
$2^4$&2  (sec) & 0.0186698 &3.4852   \\ \hline
$2^5$&3  (sec) & 0.00580343 &3.21703 \\ \hline
$2^6$&7  (sec) & 0.00187212 &3.09993\\ \hline
$2^7$&13  (sec) & 0.00060713 &3.08355 \\ \hline
$2^8$&27  (sec) & 0.000201581 &3.01185\\ \hline
\end{tabular}
\end{center}
 $$\mbox{ Table 2: Results corresponding to Example 5.2: sphere centered
   at (0,0,2) with radius 1.}$$ 
 $$\mbox{ Parameters: } F=2^7-1,\ N_d=2,\ 0\le \rho \le 4.$$

\begin{center}
\renewcommand{\arraystretch}{1.25}
\begin{tabular}{|r|p{1.6in} |r  |r|}
\hline
\textbf{$N_i$} &\textbf{time per iteration} &\textbf{error}&\textbf{error ratio}\\ \hline
$2^3$&2  (sec)  &0.0029432  & \\ \hline
$2^4$&3  (sec) & 0.00127203 & 2.31378   \\ \hline
$2^5$&7  (sec) & 0.000600751 & 2.11741 \\ \hline
$2^6$&14  (sec) & 0.000250397 &2.39919\\ \hline
$2^7$&27  (sec) & 8.41992e-05 &2.97387 \\ \hline
$2^8$&54  (sec) &1.89269e-05 & 4.44865 \\ \hline
\end{tabular}
\end{center}
$$\mbox{ Table 3: Same parameters with table 2 except } N_d=4.$$  

\begin{center}
\renewcommand{\arraystretch}{1.25}
\begin{tabular}{|r|p{1.6in} |r  |r|}
\hline
\textbf{$N_i$}&\textbf{time per iteration}&\textbf{error}&\textbf{error ratio}\\ \hline
$2^3$ & 4  (sec)&0.000757108 & \\ \hline
$2^4$ & 7  (sec)&0.000327485 &2.31189  \\ \hline
$2^5$ & 14  (sec)&0.000131013 &2.49964 \\ \hline
$2^6$ & 27  (sec)&3.06305e-05 &4.27719\\ \hline
$2^7$ & 53  (sec)&7.49014e-06 &4.08945\\ \hline
\end{tabular}
\end{center}
 $$\mbox{ Table 4: Same parameters with table 2 except } N_d=8.$$ 

\begin{center}
\renewcommand{\arraystretch}{1.25}
\begin{tabular}{|r|p{1.6in} |r  |r|}  
\hline
 \textbf{F}  &\textbf{time per iteration}&\textbf{error}& $\log_2$(\textbf{error ratio})  \\ \hline
$2^2$-1&1  (sec) &0.313367   & \\ \hline 
$2^3$-1&1  (sec) &0.0378805   &3.04833\\ \hline
$2^4$-1&2  (sec) &0.00503166  &2.91235 \\ \hline
$2^5$-1&5  (sec) &0.00100595  &2.32247  \\ \hline
$2^6$-1&11  (sec) &0.000180542 &2.47816   \\ \hline
\end{tabular}
\end{center}
$$\mbox{ Table 5: Angular convergence}$$
 $$\mbox{ Parameters: } N_i=2^{8},\ N_d=2,\ 0\le \rho \le 4.$$\\
 
\begin{example}
To test the convergence of the angular integration, our third example relates to the body of the revolution generated by an annulus. 
This hollowed sphere has the refractive index
\begin{equation}\label{25}
n(\rho,\cos\theta)=
\begin{cases}
(1+\sqrt{\mid \cos\theta \mid ^p})^{1/2}, &1\le \rho \le2\\
1, &0\le \rho <2 \mbox{ or } \rho >2.
\end{cases}
\end{equation}
Then,
$$m(\rho,t)=
\begin{cases}
-|t|^{\beta}, &1\le \rho \le2\\
0, &0\le \rho <2 \mbox{ or } \rho >2.
\end{cases}
$$

The Legendre coefficients of $m(\rho,t)\mbox{,  }1\le \rho \le2$ 
is given by $\sum_{n=0}^{\infty}m_{2n}(\rho)P_{2n}(t)$,\\ 
where
\begin{equation}\label{26}
\begin{split}
m_{2n}(\rho)&=-\frac{1}{2n+1}\int\limits_{-1}^{1}\!P_{2n}(t)|t|^{\beta}\,dt\\
&=-\frac{2}{2n+1}\int\limits_{0}^{1}\!P_{2n}(t)t^{\beta}\,dt\\
&=-\frac{(4n+1)\sqrt{\pi}2^{-\beta-1}\Gamma(1+\beta)}{\Gamma(1-n+\frac{1}{2}\beta)
\Gamma(\frac{1}{2}\beta+\frac{3}{2}+n)}\\
&=-(4n+1)\sqrt{\pi}2^{-\beta-1}\frac{\Gamma(1+\beta)}{\Gamma(1+\frac{1}{2}\beta)
\Gamma(\frac{3}{2}+\frac{1}{2}\beta)}\frac{\prod_{k=0}^{n-1}(\frac{1}{2}\beta-k)}
{\prod_{k=0}^{n-1}(\frac{1}{2}\beta+\frac{3}{2}+k)}.
\end{split}
\end{equation}  

In tables 6-9, we present the order of convergence for $\beta=$.2, 1.7, 2.2 and
 3.2 respectively. The tables clearly show 
 the correlation between smoothness of the refractive index $n(x)$ and the order of convergence.
\end{example}

\begin{center}
\renewcommand{\arraystretch}{1.25}
\begin{tabular}{|r|p{1.6in} | p{1in}    |r|}  
\hline
 \textbf{F}  & \textbf{time per iteration}&\textbf{error} & $\log_2$(\textbf{error ratio})  \\ \hline
$2^4$-1 & 1  (sec)&0.000543237 & 5.22994    \\ \hline
$2^5$-1 & 2  (sec)&0.000123797 & 2.13361  \\ \hline
$2^6$-1 & 5  (sec)&2.75278e-05 & 2.16901  \\ \hline
$2^7$-1 & 13  (sec)&5.78823e-06 & 2.2497 \\ \hline
$2^8$-1 & 30  (sec)&1.06212e-06 & 2.44617 \\ \hline
$2^9$-1 & 67  (sec)&1.96433e-07 & 2.43484 \\ \hline 
\end{tabular}
\end{center}
$$\mbox{ Table 6: }\beta=.2,\ m\in C^{.2},\ u\in C^{2.2}$$
$$\mbox{ Parameters: } F=2^{10}-1,\ N_d=2,\ N_i=2^7,\ 0\le \rho \le 4.$$\\

\begin{center}
\renewcommand{\arraystretch}{1.25}
\begin{tabular}{|r| p{1in}    |r|}
\hline
 \textbf{F}  & \textbf{error}& $\log_2$(\textbf{error ratio})   \\ \hline
$2^4$-1 &1.70623e-05 &    \\ \hline
$2^5$-1 & 1.10166e-06 &3.95306  \\ \hline
$2^6$-1 & 8.65346e-08 &3.67025  \\ \hline
$2^7$-1 & 6.6141e-09  & 3.70966 \\ \hline
$2^8$-1 & 4.5934e-10 & 3.84791 \\ \hline
$2^9$-1 & 3.25004e-11 & 3.82104 \\ \hline 
\end{tabular}
\end{center}
$$\mbox{ Table 7: }\beta=1.7,\ m\in C^{1.7},\ u\in C^{3.7}$$
$$\mbox{ Same parameters with table 6.}$$\\ 

\begin{center}
\renewcommand{\arraystretch}{1.25}
\begin{tabular}{|r| p{1in}  |r|}
\hline
 \textbf{F}  & \textbf{error} & $\log_2$(\textbf{error ratio}) \\ \hline
$2^4$-1 & 4.85002e-06 & \\ \hline
$2^5$-1 & 1.87781e-07 & 4.69086 \\ \hline
$2^6$-1 & 1.03671e-08 & 4.17896 \\ \hline
$2^7$-1 & 5.61064e-10 & 4.20771 \\ \hline
$2^8$-1 & 2.78307e-11 & 4.33342 \\ \hline
$2^9$-1 & 1.40967e-12 & 4.30325 \\ \hline 
\end{tabular}
\end{center}
$$\mbox{ Table 8: }\beta=2.2,\ m\in C^{2.2 },\ u\in C^{4.2}$$
$$\mbox{ Same parameters with table 6.}$$ \\

\begin{center}
\renewcommand{\arraystretch}{1.25}
\begin{tabular}{|r| p{1in}  |r|}
\hline
 \textbf{F}  & \textbf{error} & $\log_2$(\textbf{error ratio}) \\ \hline
$2^4$-1 &  3.43503e-06 & \\ \hline
$2^5$-1 &  5.82177e-08 & 5.88272\\ \hline
$2^6$-1 &  1.32625e-09 & 5.45603\\ \hline
$2^7$-1 &  3.58185e-11 & 5.21051\\ \hline
$2^8$-1 &  8.9928e-13 & 5.31579\\ \hline
$2^9$-1 &  2.3999e-14 & 5.22773\\ \hline 
\end{tabular}
\end{center}
$$\mbox{ Table 8: }\beta=3.2,\ m\in C^{3.2 },\ u\in C^{5.2}$$
$$\mbox{ Same parameters with table 6.}$$ \\

\section{Numerical analysis}
Actually, what we have solved is following,
\begin{equation}\label{27}
v(\rho,\cos \theta)=u^{i,F}(\rho,\cos \theta)+\frac{i}{2}(K^Fv)(\rho,\cos \theta),
\end{equation}
where
\begin{equation}\label{28}
u^{i,F}(\rho,\cos \theta)=\sum^F_{n=0}u^i_n(\rho)P_n(\cos\theta)
\end{equation}
\begin{equation}\label{29}
(K^Fv)(\rho,\cos \theta)=\sum^F_{n=0}(K_nv)(\rho)P_n(\cos\theta)
\end{equation}
Therefore,
$$
  v_n(\rho)=
\begin{cases}
u_n^i(\rho)+\frac{i}{2}(K_nv)(\rho),&n \le F\\ 
0,&n > F.\\
\end{cases}
$$

While the exact solution satisfies,
\begin{equation}\label{30}
u^F(\rho,\cos \theta)=u^{i,F}(\rho,\cos \theta)+\frac{i}{2}(K^Fu)(\rho,\cos \theta),
\end{equation}
Therefore, the error in the solution of the approximate integral
equation (\ref{27}) is given by
\begin{equation}\label{31}
|u(\rho,\cos\theta)-v^F(\rho,\cos\theta)| \le |u^F(\rho,\cos\theta)-v^F(\rho,\cos\theta)|+|u^T(\rho,\cos\theta)|, 
\end{equation}
where 
$$
u^T(\rho,\cos\theta)=u(\rho,\cos\theta)-u^F(\rho,\cos\theta)=\sum_{n>F}u_n(\rho)P_n(\cos\theta)
$$
By comparing (\ref{27}) and(\ref{30}), we get $$
  u^F-v^F=\frac{i}{2}K^F(u^F-v^F)+\frac{i}{2}K^Fu^T
$$
Therefore,
\begin{equation*}
\begin{split}
\|u^F-v^F\|_{\infty}
&=\|(I-\frac{i}{2}K^F)^{-1}\|_{\infty}\|\frac{i}{2}K^Fu^T\|_{\infty}\\
&\le B\|K^Fu^T\|_{\infty}\\
\end{split}                    
\end{equation*}
where,
$$\|(I-\frac{i}{2}K^F)^{-1}\|_{\infty}\le B.$$
From the equations (\ref{13}), (\ref{14}),
\begin{equation}\label{32}
 (K_{n}u^T)(\rho,\cos\theta)=-(2n+1)k^3\int\limits_{0}^{\mathbf R}\int\limits_{0}^{\pi}\!
h_{n}^{(1)}(k\rho_{>}) j_{n}(k\rho_{<})\rho^2(u^{F}(\rho,\theta) m ^{2F}(\rho,\theta)P_{n}(\cos\theta)\sin\theta\,d\theta)\,d\rho,
\end{equation}

\begin{lemma}\label{lm6.2}
There is a constant $C(R,k)$ such that
\begin{equation*}
\| \int\limits_{0}^{R}\!h_{n}^{(1)}(k\rho_{>})
j_{n}(k\rho_{<})\rho^2 \,d\rho \|_{\infty}\le \frac{C(R,k)}{max(1,n^2)}  
\end{equation*}
\end{lemma}

\begin{proof}
\begin{equation*}
\begin{split}
\| \int\limits_{0}^{R}\!h_{n}^{(1)}(k\rho_{>})
j_{n}(k\rho_{<})\rho^2 \,d\rho \|_{\infty}
 &\le |y_{n}(k\widetilde{\rho})|\int\limits_{0}^{\min(\widetilde{\rho},\mathbf R)}\!|j_{n}(k\rho)|
  \rho^2\,d\rho
  +|j_{n}(k\widetilde{\rho})|\int\limits_{\min(\widetilde{\rho},\mathbf R)}^{\mathbf R}\!|y_{n}(k\rho)| 
  \rho^2\,d\rho  \\
 & +|j_{n}(k\widetilde{\rho})|\int\limits_{0}^{\mathbf
   R}\!|j_{n}(k\rho)|\rho^2\,d\rho. \\
\end{split}
\end{equation*}
\begin{equation*}
\begin{split}
&j_n(\rho)=\frac{\rho^n}{1\cdot3\cdot5\cdot\dots(2n+1)}
\bigg[1-\frac{\frac{1}{2}\rho^2}{1!(2n+3)}+\frac{({\frac{1}{2}\rho^2})^2}{2!(2n+3)(2n+5)}+\dots
\bigg]\\
&y_n(\rho)=\frac{-1\cdot1\cdot3\cdot5\cdot\dots(2n-1)}{\rho^{n+1}}
\bigg[1-\frac{\frac{1}{2}\rho^2}{1!(1-2n)}+\frac{({\frac{1}{2}\rho^2})^2}{2!(1-2n)(3-2n)}+\dots
\bigg]\\
\end{split}
\end{equation*}
Therefore,
\begin{equation*}
|j_n(ka)y_n(kb)|\le\frac{a^n}{b^{n+1}(2n+1)}C(R,k)
\end{equation*}
Then, the result follows immediately.
\end{proof}\\

Therefore,
\begin{equation}\label{33}
\begin{split}
\|K_{n}u^T\|_{\infty}& \le (2n+1)k^3(\int\limits_{0}^{\mathbf
  R}\!|h_{n}^{(1)}(k\rho_{>}) j_{n}(k\rho_{<})|\rho^2
\,d\rho)\sum_{s>F}\sum_{l=s-n}^{s+n}
\|m_l(\rho)\|_{\infty}\|u_s(\rho)\|_{\infty}\int\limits_{-1}^{1}\!|P_mP_sP_n|\,d\rho\\
&\le \frac{C(k,R)}{max(1,n)}\sum_{s>F}\sum_{l=s-n}^{s+n}
\|m_l(\rho)\|_{\infty}\|u_s(\rho)\|_{\infty}\int\limits_{-1}^{1}\!|P_mP_sP_n|\,d\rho\\
\end{split}
\end{equation}
where,
\begin{equation}\label{34}
\begin{split}
  &m(\rho,\cos\theta)=\sum_{m=0}^{\infty}m_l(\rho)P_l(\cos\theta),\\
  &u^T(\rho,\cos\theta)=\sum_{s>F}u_s(\rho)P_s(\cos\theta),\\
\end{split}
\end{equation}
\begin{lemma}
If $f\in C([-1,1])$, the the Legendre coefficients of $f$,
\begin{equation}\label{35}
c_n=\frac{2n+1}{2}\int\limits_{-1}^{1}\!f(t)P_n(t)\,dx \\
= \frac{1}{n+1}\frac{2n+1}{2}\int\limits_{-1}^{1}\frac{df(t)}{dt}
(P_{n-1}(t)-tP_n(t)) \!\,dt.\\
\end{equation}
\end{lemma}

\begin{proof}
Since, $\varTheta=P_n(\cos\theta)$ satisfies, 
\begin{equation*}
  \frac{d}{d\theta}\big(\sin\theta\frac{d\varTheta}{d\theta}\big)=-n(n+1)(\sin\theta)\varTheta,
\end{equation*}
Therefore, 
\begin{equation*}
\begin{split}
\frac{2}{2n+1}c_n&=\int\limits_{-1}^{1}\!f(x)P_n(x)\,dx
=\int\limits_{0}^{\pi}\!f(\cos\theta)P_n(\cos\theta)\sin\theta\,d\theta\\
&=-\int\limits_{0}^{\pi}\!\frac{f(\cos\theta)}{n(n+1)}\frac{d}{d\theta}\big(\sin\theta
\frac{d(P_n(\cos\theta))}{d\theta}\big)d\theta\\
&=-\frac{f(\cos\theta)}{n(n+1)}\sin\theta
\frac{d(P_n(\cos\theta))}{d\theta}\Big|^{\pi}_{0}+\frac{1}{n(n+1)}
\int\limits_{0}^{\pi}\frac{d}{d\theta}(f(\cos\theta))\cdot\sin\theta
\frac{d(P_n(\cos\theta))}{d\theta} \!\,d\theta\\
&=\frac{1}{n(n+1)}
\int\limits_{0}^{\pi}\frac{d}{d\theta}(f(\cos\theta))\cdot\sin\theta
\frac{d(P_n(\cos\theta))}{d\theta} \!\,d\theta\\
&=\frac{1}{n(n+1)}\int\limits_{-1}^{1}\frac{df(t)}{dt}\frac{dt}{d\theta}
\frac{dP_n(t)}{dt}\frac{dt}{d\theta} \!\,dt \\
&=\frac{1}{n(n+1)}\int\limits_{-1}^{1}\frac{df(t)}{dt}
\frac{dP_n(t)}{dt}(1-t^2) \!\,dt\\
&= \frac{1}{n+1}\int\limits_{-1}^{1}\frac{df(t)}{dt}
(P_{n-1}(t)-tP_n(t)) \!\,dt\\
\end{split}
\end{equation*}
For the last equality, we used the formula, $\frac{dP_n(t)}{dt}(1-t^2)=n(P_{n-1}(t)-tP_n(t))$.
\end{proof}

Stieltjes' formula (\cite{sansone}) says,
\begin{equation}\label{36}
\begin{split}
P_n(\cos\theta)&=\frac{4}{\pi}\frac{1}{(2n+1)}\frac{1}{\alpha_n}\bigg[
\frac{\cos\big((n+\frac{1}{2})\theta-\frac{\pi}{4}\big)}{(2\sin\theta)^{\frac{1}{2}}}\\
&+\alpha_1\frac{1}{2n+3}\frac{\cos\big((n+\frac{3}{2})\theta-\frac{3\pi}{4}\big)}{(2\sin\theta)^{\frac{3}{2}}}\\
&\alpha_2\frac{1\cdot3}{(2n+3)(2n+5)}\frac{\cos\big((n+\frac{5}{2})\theta-\frac{5\pi}{4}\big)}{(2\sin\theta)^{\frac{5}{2}}}+\cdots \bigg].\\
\end{split}
\end{equation}
where
$$
 \alpha_n=\frac{(2n-1)!!}{(2n)!!} 
$$

The series is convergent for $\frac{\pi}{6}<\theta<\frac{5\pi}{6}$. If
the sum of the first $m$ terms is taken as an approximation of
$P_n(\cos\theta)$ for $0<\theta<\pi$, the error $p_{n,m}(\theta)$
satisfies the inequality
\begin{equation}\label{37}
\begin{split}
|p_{n,m}|<&\frac{4}{\pi}\frac{1}{(2n+1)}\frac{1}{\alpha_n}\\
&\cdot \alpha_m\frac{(2m-1)!!}{(2n+3)(2n+5)\dots(2n+2m+1)}
\frac{2}{(2\sin\theta)^{m+\frac{1}{2}}}.\\
\end{split}
\end{equation}
In other words, the error is less than twice the absolute value of the
first term omitted if the cosine term is replaced by 1.
\begin{lemma}
If $\frac{df}{dt}$ is integrable, then the Legendre coefficients of $f$,
\begin{equation}\label{38}
|c_n|=|\frac{2n+1}{2}\int\limits_{-1}^{1}\!f(x)P_n(x)|\,dx\le \frac{C}{n^{\frac{1}{2}}}\\
\end{equation}
\end{lemma}
\begin{proof}
Since $P_n(\cos\theta)$ satisfies, 
\begin{equation*} 
\Bigg|P_n(\cos\theta)-\frac{4}{\pi}\frac{1}{(2n+1)}\frac{1}{\alpha_n}\bigg[
\frac{\cos\big((n+\frac{1}{2})\theta-\frac{\pi}{4}\big)}{(2\sin\theta)^{\frac{1}{2}}}\bigg]\Bigg|
\le \frac{4}{\pi}\frac{1}{(2n+1)(2n+3)}\frac{1}{\alpha_n}\frac{1}{(2\sin\theta)^{\frac{3}{2}}}
\end{equation*}

and, 
\begin{equation*}
\begin{split}
\int\limits_{0}^{\pi}\!f(\cos\theta)\cos(n\theta)(\sin\theta)^{\frac{1}{2}}\,d\theta&=
f(\cos\theta)(\sin\theta)^{\frac{1}{2}}\frac{\sin(n\theta)}{n}\Big|_0^{\pi}-
\frac{1}{n}\int\limits_{0}^{\pi}\!
\frac{d\big(f(\cos\theta)(\sin\theta)^{\frac{1}{2}}\big)}{d\theta}\sin(n\theta) \,d\theta\\
&=O\big(\frac{1}{n}\big)\\
\end{split}
\end{equation*}
Similarly,
\begin{equation*}
\int\limits_{0}^{\pi}\!f(\cos\theta)\sin(n\theta)(\sin\theta)^{\frac{1}{2}}\,d\theta=O\big(\frac{1}{n}\big)
\end{equation*}
Therefore,
\begin{equation*}
\begin{split}
|c_n|&\frac{2}{2n+1}=|\int\limits_{-1}^{1}\!f(t)P_n(t)|\,dt\\
&\le \Big|\int\limits_{0}^{\pi}\!f(\cos\theta)\frac{4}{\pi}\frac{1}{(2n+1)}\frac{1}{\alpha_n}
\cos\big((n+\frac{1}{2})\theta-\frac{\pi}{4}\big)\frac{(\sin\theta)^{\frac{1}{2}}}{\sqrt2}\,d\theta
\Big|+\frac{4}{\pi}\frac{1}{(2n+1)(2n+3)}\frac{1}{\alpha_n}
\int\limits_{0}^{\pi}\!\frac{\,d\theta}{2^\frac{2}{3}(\sin\theta)^\frac{1}{2}} \\
&\le O\big(\frac{1}{n^{\frac{3}{2}}}\big)
\end{split}
\end{equation*}
where
$$
 \alpha_n \thickapprox O(\frac{1}{\sqrt n})
$$
\end{proof}
\begin{corollary}
If $\frac{d^mf}{dt^m}$ is integrable, then the Legendre coefficients of $f$,
\begin{equation}\label{39}
|c_n|=|\frac{2n+1}{2}\int\limits_{-1}^{1}\!f(t)P_n(t)|\,dx\le \frac{C}{n^{m-\frac{1}{2}}}\\
\end{equation}
\end{corollary}
\begin{proof}
Combining (\ref{35}) and (\ref{38}), it follows directly.
\end{proof}\\

For example, if $f(t)=|t|^{k+\alpha}\in \mbox{  } C^{k+\alpha}$, $0< \alpha <1$, then the
Legendre coefficients decay like $O(\frac{1}{n^{k+.5}})$ by the corollary. But, numerically, it shows that  
$f(t)$ decays as $ O(\frac{1}{n^{k+.5+\alpha}})$. But,
I will not attempt to prove the correlation between $\alpha$ and the
convergence order here.


\begin{lemma}
Suppose $\{m_l(\rho)\},\mbox{ }\{u_s(\rho)\}$ from the equations (\ref{34}) decay like 
$O(\frac{1}{l^{k+2+\alpha}}),\mbox{ }O(\frac{1}{s^{k+4+\alpha}})$ with $\alpha > 0$ respectively. 
If $k\ge 1$ then,
 $$\|K^{F}u^T\|_{\infty}=O(\frac{1}{F^{k+4.5+\alpha+min(1.5,k+\alpha)}}).$$ 
If $k=0$ then,
$$\|K^{F}u^T\|_{\infty}=O(\frac{\log F}{F^{4.5+\alpha+min(1,.5+\alpha)}}).$$ 
\end{lemma}
\begin{proof}
If $k\ge 1$,
\begin{equation}\label{40}
\begin{split}
\|K_{n}u^T\|_{\infty}
&\le \frac{C(k,R)}{max(1,n)}\sum_{s>F}\sum_{l=s-n}^{s+n}
\|m_l(\rho)\|_{\infty}\|u_s(\rho)\|_{\infty}\int\limits_{-1}^{1}\!|P_mP_sP_n|\,d\rho\\
&\le \frac{C(k,R)}{max(1,n)}\sum_{s>F}\sum_{l=s-n}^{s+n}\frac{1}{l^{k+2+\alpha}}
\frac{1}{s^{k+4+\alpha}}\int\limits_{-1}^{1}\!|P_mP_sP_n|\,d\rho\\
&\le
\frac{C(k,R)}{max(1,n)}\sum_{s>F}\frac{1}{(s-n)^{k+1+\alpha}}\frac{1}{s^{k+4.5+\alpha}}
\frac{1}{max(1,n^{.5})}\\
&\le  \frac{C(k,R)}{max(1,n^{1.5})}\frac{1}{(F+1-n)^{k+\alpha}}\frac{1}{F^{k+4.5+\alpha}}\\
\end{split}
\end{equation}
Therefore,
\begin{equation}\label{41}
\begin{split}
\|K^{F}u^T\|_{\infty}&\le \sum_{n=0}^{F}\|K_{n}u^T\|_{\infty}\\
&\le \frac{C(k,R)}{F^{k+4.5+\alpha}}\sum_{n=0}^{F}
\frac{1}{max(1,n^{1.5})}\frac{1}{(F+1-n)^{k+\alpha}}
\end{split}
\end{equation}
Since 
\begin{equation}\label{42}
\begin{split}
\sum_{n=1}^{F}
\frac{1}{n^{1.5}}&\frac{1}{(F+1-n)^{k+\alpha}}\\
&\le \sum_{n=1}^{\ulcorner \frac{F}{2} \urcorner}\frac{1}{n^{1.5}(F+1-n)^{k+\alpha}}+
\sum_{n=1}^{\ulcorner \frac{F}{2} \urcorner}\frac{1}{n^{k+\alpha}(F+1-n)^{1.5}}\\ 
&=O(\frac{1}{F^{k+\alpha}})+O(\frac{1}{F^{min(k+\alpha+.5,1.5)}})
\le O(\frac{1}{F^{min(k+\alpha,1.5)}})\\
\end{split}
\end{equation}
Combining equations (\ref{41}) and (\ref{42}) give our desired result. \\
If $k=0$,
\begin{equation}\label{43}
\begin{split}
\|K_{n}u^T\|_{\infty}
&\le \frac{C(k,R)}{max(1,n)}\sum_{s>F}\sum_{l=s-n}^{s+n}\frac{1}{l^{2+\alpha}}
\frac{1}{s^{4+\alpha}}\int\limits_{-1}^{1}\!|P_mP_sP_n|\,d\rho\\
&\le \frac{C(k,R)}{max(1,n)}\sum_{s>F}\frac{1}{(s-n)^{1.5+\alpha}}\frac{1}{s^{4.5+\alpha}}\\
&\le  \frac{C(k,R)}{max(1,n)}\frac{1}{(F+1-n)^{.5+\alpha}}\frac{1}{F^{4.5+\alpha}}\\
\end{split}
\end{equation}
Therefore,
\begin{equation}\label{44}
\begin{split}
\|K^{F}u^T\|_{\infty}
&\le \frac{C(k,R)}{F^{4.5+\alpha}}\sum_{n=0}^{F}
\frac{1}{max(1,n)}\frac{1}{(F+1-n)^{.5+\alpha}}\\
&\le \frac{C(k,R)}{F^{4.5+\alpha}}(\sum_{n=1}^{\ulcorner \frac{F}{2} \urcorner}\frac{1}{n(F+1-n)^{.5+\alpha}}+
\sum_{n=1}^{\ulcorner \frac{F}{2} \urcorner}\frac{1}{n^{.5+\alpha}(F+1-n)}+O(\frac{1}{F^{.5+\alpha}}))\\
&\le \frac{C(k,R)}{F^{4.5+\alpha}}(O(\frac{\log F}{F^{.5+\alpha}})+
O(\frac{\log F}{F^{min(1,.5+\alpha)}})+O(\frac{1}{F^{.5+\alpha}}))
=O(\frac{\log F}{F^{4.5+\alpha+min(1,.5+\alpha)}}) 
\end{split}
\end{equation}

\end{proof}
\begin{lemma}
If $\{m_l(\rho)\},\mbox{ }\{u_s(\rho)\}$ decay like 
$O(\frac{1}{l^{k+\alpha}}),\mbox{ }O(\frac{1}{s^{k+2+\alpha}})$ respectively
where $k=0,\mbox{ }1$. 
Then, 
$$
\|K^{F}u^T\|_{\infty}=
\begin{cases}
O(\frac{1}{F^{3+2\alpha}}) & \mbox{ if  }k=1,\\
O(\frac{1}{F^{1+2\alpha}}) & \mbox{ if  }k=0.\\
\end{cases}
$$
\end{lemma}

\begin{proof}
If $k=1$,
\begin{equation}\label{45}
\begin{split}
\|K_{n}u^T\|_{\infty}
&\le \frac{C(k,R)}{max(1,n)}\sum_{s>F}\sum_{l=s-n}^{s+n}\frac{1}{l^{1+\alpha}}
\frac{1}{s^{3+\alpha}}\int\limits_{-1}^{1}\!|P_mP_sP_n|\,d\rho\\
&\le C(k,R)\sum_{s>F}\frac{1}{(s-n)^{1+\alpha}}\frac{1}{s^{3.5+\alpha}}\frac{1}{max(1,n^{.5})}\\
&\le  \frac{C(k,R)}{max(1,n^{.5})}\frac{1}{(F+1-n)^{\alpha}}\frac{1}{F^{3.5+\alpha}}\\
\end{split}
\end{equation}
Similarly,
\begin{equation}\label{46}
\begin{split}
\sum_{n=1}^{F}
\frac{1}{n^{.5}}&\frac{1}{(F+1-n)^{\alpha}}\\
&\le \sum_{n=1}^{\ulcorner \frac{F}{2} \urcorner}\frac{1}{n^{.5}(F+1-n)^{\alpha}}+
\sum_{n=1}^{\ulcorner \frac{F}{2} \urcorner}\frac{1}{n^{\alpha}(F+1-n)^{.5}}\\ 
&=O(\frac{1}{F^{-.5+\alpha}})+O(\frac{1}{F^{min(-.5+\alpha,.5)}})
\le O(\frac{1}{F^{-.5+\alpha}})\\
\end{split}
\end{equation}
Therefore,
\begin{equation}\label{47}
\|K^{F}u^T\|_{\infty}=O(\frac{1}{F^{3+2\alpha}})
\end{equation}
If $k=0$, then
\begin{equation}\label{48}
\begin{split}
\|K_{n}u^T\|_{\infty}
&\le \frac{C(k,R)}{max(1,n)}\sum_{s>F}\sum_{l=s-n}^{s+n}\frac{1}{l^{\alpha}}
\frac{1}{s^{2+\alpha}}\int\limits_{-1}^{1}\!|P_mP_sP_n|\,d\rho\\
&\le C(k,R)\sum_{s>F}\frac{1}{(s-n)^{\alpha}}\frac{1}{s^{2.5+\alpha}}\frac{1}{max(1,n^{.5})}\\
&\le  \frac{C(k,R)}{max(1,n^{.5})}\frac{1}{(F+1-n)^{\alpha}}\frac{1}{F^{1.5+\alpha}}\\
\end{split}
\end{equation}
Therefore,
\begin{equation}\label{49}
\|K^{F}u^T\|_{\infty}=O(\frac{1}{F^{1+2\alpha}})
\end{equation}
\end{proof}

\begin{corollary}
If $\{m_l(\rho)\},\mbox{ }\{u_s(\rho)\}$ decay like 
$O(\frac{1}{m^{k+\alpha}}),\mbox{ }O(\frac{1}{s^{k+2+\alpha}})$ respectively. 
Then, 
$$
\varepsilon_F:=\|u_F-v_F\|_{\infty}=
\begin{cases}
O(\frac{1}{F^{k+2.5+\alpha+min(1.5,k-2+\alpha)}}) & \mbox{ if  }k\ge 3\\
O(\frac{\log F}{F^{4.5+\alpha+min(1,.5+\alpha)}}) & \mbox{ if  }k=2,\\
O(\frac{1}{F^{3+2\alpha}}) & \mbox{ if  }k=1,\\
O(\frac{1}{F^{1+2\alpha}}) & \mbox{ if  }k=0.\\
\end{cases}
$$ 
\end{corollary}
From (\ref{31}), the total error is given as
\begin{equation}\label{50}
\|u(\rho,\cos\theta)-v^F(\rho,\cos\theta)\|_{\infty}
 \le \|u^F(\rho,\cos\theta)-v^F(\rho,\cos\theta)\|_{\infty}+\|u^T(\rho,\cos\theta)\|_{\infty}, 
\end{equation}
where 
$$
u^T(\rho,\cos\theta)=u(\rho,\cos\theta)-u^F(\rho,\cos\theta)=\sum_{n>F}u_n(\rho)P_n(\cos\theta).
$$
So far, we have error bound for $\|u_F-v_F\|_{\infty}$ and the error bound
for $\|u^T(\rho,\cos\theta)\|_{\infty}$ is easily given as
$$
\varsigma_F:=\|\sum_{n>F}u_n(\rho)P_n(\cos\theta)\|_{\infty}=
O(\frac{1}{F^{k+1+\alpha}}) 
$$
Where $u_n(\rho)\thicksim O(\frac{1}{n^{k+2+\alpha}})$ 
Therefore, the total error is dominated by $\varsigma_F$.\\ 

If $m(\rho,\theta)$ is given as (\ref{25}) then, $m_n(\rho) \thicksim \frac{1}{n^{.5+\beta}}$ and $u_n(\rho)\thicksim O(\frac{1}{n^{2.5+\beta}})$. 
Therefore, the total error is given as
$$\|u(\rho,\cos\theta)-v^F(\rho,\cos\theta)\|_{\infty}\thicksim O(\frac{1}{F^{1.5+\beta}})$$.
But, actually, our tables 6-9 show $\|u-v^F\|_{\infty}\thicksim O(\frac{1}{F^{2+\beta}})$. 

The reason is that 
$\varsigma_F=\|\sum_{n>F}u_n(\rho)P_n(\cos\theta)\|_{\infty} $ 
converges faster than $\varphi_F:=\|\sum_{n>F}|u_n(\rho)|\|_{\infty}$ and the latter is $O(\frac{1}{F^{1.5+\beta}})$. For example, when $\beta=2.2$,
\begin{center}
\renewcommand{\arraystretch}{1.25}
\begin{tabular}{|r|r| p{1in}    |r| p{1in}    |r|}  
\hline
 \textbf{F}  & \textbf{total error} & $\varsigma_F$  & $\log_2$($\varsigma_F$) 
& $\varphi_F$ & $\log_2$($\varphi_F$)   \\ \hline
$2^4$-1 & 4.85002e-06   & 4.82746e-06 &        & 2.4488e-05 &   \\ \hline
$2^5$-1 &  1.87781e-07& 1.87812e-07 &4.6839  &1.72411e-06 &3.82815\\ \hline
$2^6$-1 & 1.03671e-08 & 1.03674e-08 & 4.17917 &1.28051e-07 & 3.75106\\ \hline
$2^7$-1 & 5.61064e-10 & 5.61065e-10 & 4.20774 &9.72137e-09  &3.71942\\ \hline
$2^8$-1 & 2.78307e-11 & 2.78301e-11 & 4.33345 &7.44753e-10 & 3.70633\\ \hline
$2^9$-1 & 1.40967e-12 & 1.40934e-12 & 4.30355 &5.29085e-11 &3.81519\\ \hline 
\end{tabular}
\end{center}
$$\mbox{ Table 10: }\beta=2.2,\ m\in C^{2.2},\ u\in C^{4.2}$$
$$\mbox{ Parameters: } F=2^{10}-1,\ N_d=2,\ N_i=2^7,\ 0\le \rho \le 4.$$\\
 


\appendix
\section{Appendix}
\subsection{The Fast Legendre Transform and its Inverse}

\subsubsection{General theory of orthogonal polynomials}
\begin{lemma}(three-term recurrence) Let  $\{p_{m}\}$ be an orthogonal polynomial sequence for a 
nonnegative integrable weight function. Then $\{p_{m}\}$ satisfies a three-term recurrence relation
\begin{equation}\label{55}
p_{k+1}(x)=(A_{k}x+B_{k})p_{k}(x)+C_{k}p_{k-1}(x),
\end{equation}
where $A_k$, $B_k$, $C_k$ are real numbers with $A_k\ne0$ and $B_k\ne0$.
\end{lemma}
\begin{proof}
See, e.g., [\cite{gordon}, Theorem 4.1]
\end{proof}

Next define the $associated$  $polynomials$ $Q_{l,m},\mbox{  }R_{l,m}$ for the orthogonal polynomial sequence $\{p_m\}$ by the following recurrences on m \cite{associated}
\begin{equation}\label{56}
\begin{split}
Q_{l,m}(x) =& (A_{l+m-1}x+B_{l+m-1})Q_{l,m-1}(x)+C_{l+m-1}Q_{l,m-2}(x),\\
            & Q_{l,0}(x)=1,\mbox{ }Q_{l,1}=A_lx+B_l, \\
R_{l,m}(x) =&(A_{l+m-1}x+B_{l+m-1})R_{l,m-1}(x)+C_{l+m-1}R_{l,m-2}(x),\\
            & R_{l,0}(x)=0,\mbox{ }R_{l,1}=C_l.
\end{split}
\end{equation}

We have
\begin{lemma}\label{lm8.1}(generalized three-term recurrence).
The associated polynomials satisfy 
$\deg Q_{l,m}=m$, $\deg R_{l,m}\leq m-1,$ and for $l\ge 1$ 
and $m \ge 1,$
\begin{equation}\label{57}
p_{l+m}=Q_{l,m}\cdot p_l+R_{l,m}\cdot p_{l-1}.
\end{equation}
\end{lemma}
\begin{proof}
Equation (\ref{57}) follows by induction on $m$ with the case $m=1$ being 
the original three-term recurrence (\ref{55}).
\end{proof}


\subsubsection{The Legendre Transform} 
Assume that $f(x)$ is a polynomial of degree less than $N$ and let
$$ 
f(x)=\sum_{n=0}^{N-1} c_n T_n(x).
$$
where
\begin{equation}\label{52}
\begin{split}
  c_{n}= & \frac{(-1)^{n}\epsilon_{n}}{N}\sum_{j=0}^{N-1}f_{j}\cos{\frac{(2j+1)n\pi}{2N}} 
\end{split}
\end{equation}
and 
$$
\epsilon_{n}=
\begin{cases}
           1,& \mbox{ if }n=0,\\
           2,& \mbox{ if }n>0. 
\end{cases}
$$ 
Then,
$$ 
\int\limits_{-1}^{1}\!f(x)\,dx
=\int\limits_{0}^{\pi}\!f(\cos\theta )\sin\theta \,d\theta 
=\sum_{n=0}^{N-1}c_n \alpha_n,
$$ 
where
$$
\alpha_{n}=\int\limits_{0}^{\pi}\! \cos(\mbox{n}\theta) \sin\theta \,d\theta.
$$
Then with (\ref{52}) we obtain,

\begin{equation}\label{59}
\begin{split}
\int\limits_{-1}^{1}\!f(x)\,dx
  & =\sum_{n=0}^{N-1}\frac{(-1)^n\epsilon_{n}}{N}\alpha_{n}\sum_{i=0}^{N-1}f(\cos\frac{(2i+1)\pi}{2N})\cos(\mbox{n}\frac{(2i+1)\pi}{2N}{i})\\
  & = \sum_{i=0}^{N-1}f(\cos\frac{(2i+1)\pi}{2N})\sum_{n=0}^{N-1}\frac{(-1)^n\epsilon_{n}}{N}\alpha_{n}\cos(\mbox{n}\frac{(2i+1)\pi}{2N}) \\
  & =\sum_{i=0}^{N-1}f(\cos \frac{(2i+1)\pi}{2N})w_{i}^{N}=\sum_{i=0}^{N-1}f(x_{i}^N)w_{i}^{N},
\end{split}
\end{equation}
where
$$
w_{i}^{N}=\sum_{n=0}^{N-1}\frac{(-1)^n\epsilon_{n}}{N}\alpha_{n}\cos(\mbox{n}\frac{(2i+1)\pi}{2N}).
$$

In particular, if the degree of $f(x)$ is less than or equal to N, then
$$
\int\limits_{-1}^{1}\!p_{m}(x)f(x)\,dx=\sum_{i=0}^{2N-1}p_{m}(x_i^{2N})f(x_i^{2N})w_{i}^{2N},
\ \mbox{m}\in {0,1,2, ,N-1}.
$$ 

\begin{definition}
The Legendre coefficients $\{c_{m}\}$ of $f(x)$ are defined as
$$
c_m=\frac{1}{\tau_{m}}\sum_{i=0}^{2N-1} P_{m}(x_i^{2N}) f(x_i^{2N}) w_{i}^{2N},
\ \mbox{m}\in {0,1,2, ,N-1},
$$ 
where $\{P_m\}$ are the Legendre polynomials and 
$\tau_{m}=\int\limits_{-1}^{1}\!{\mid P_{m}(x)\mid}^2\,dx=\frac{2}{2m+1}$. 
We define the Legendre Transform as the map
$$
  \{ f(x_i^{2N})\}_{i=0}^{2N-1} \rightarrow \{c_m\}_{m=0}^{N-1}.
$$
\end{definition}
\subsubsection{The Fast Legendre Transform (FLT): Healy-Driscoll algorithm \cite{healy}}
In \cite{flt}, Inda, Bisseling and Maslen reviewed and implemented an original idea of Healy-Driscoll
\cite{healy} to compute
$$ 
\sum_{i=0}^{N-1}p_{m}(x_i^N)f(x_i^N), \mbox{ m}\in \{0,1,2, ,N-1\}
$$ 
in $O(N\log^2 N)$ operations in the case where $\{p_m\}$ is any sequence of orthogonal polynomials, 
$x_{i}^N$ are Chebyshev points and $N$ is power of $2$.

Here we recall the ideas behind the scheme.
First we define the truncation operator $\mathcal{T}_{n}$ by 
\begin{equation}\label{60}
  \mathcal{T}_{n}f=\sum_{k=0}^{n-1}b_{k} T_{k},    
\end{equation}
where
$$
f=\sum_{k=0}^{\infty}b_{k} T_{k},  
$$
and the truncated polynomial $Z_l^K$ by
\begin{equation}\label{61}
  Z_l^K=\mathcal{T}_{K}(f\cdot p_{l}),
\end{equation}
where $\{p_{l}\}$ are orthogonal polynomials.
We also let $\mathcal{S}_{n}$ denote the 
Lagrange interpolation operator $\mathcal{S}_{n}$ that is,
$\mathcal{S}_{n}f(x)$ is the polynomial of degree less than n which agrees with 
$f(x)$ at the Chebyshev points $x_{0}^n,\dots x_{n-1}^{n}$.\\

If the degree of $f(x)$ is less than $m$, then for $n \le m$ we have 
\begin{equation}\label{62}
 \mathcal{T}_{n} f(x)=\sum_{k=0}^{n-1}b_{k} T_{k},
\end{equation}
where
$$
  \{b_{k}\}_{k=0}^{m-1}=\mbox{FCT}(\{f(x_j^m)\}_{j=0}^{m-1}).
$$
Therefore, for an arbitrary polynomial $f(x)$ and $n \le m$,
 \begin{equation}\label{63}
\mathcal{T}_{n}(\mathcal{S}_{m} f(x))=\sum_{k=0}^{n-1}b_{k} T_{k},
\end{equation} 
where 
$$
  \{b_{k}\}_{k=0}^{m-1}=\mbox{FCT}(\{f(x_j^m)\}_{j=0}^{m-1}).
$$
When the degree of the polynomial $f(x)p_l(x)$ is less than $2N$, due to the Gaussian 
quadrature (\ref{54}) we have 
\begin{equation}\label{64}
Z_l^1=\mathcal{T}_{1}(f\cdot p_{l})=\frac{1}{\pi}\int\limits_{-1}^{1}\! \frac{f(x)p_{l}dx}{\sqrt{1-x^2}} \,dx=\frac{1}{N}\sum_{i=0}^{N-1}p_{l}(x_i^N)f(x_i^N).
\end{equation}

Therefore a fast Legendre transform will follow from a scheme that computes  $Z_l^1$ fast, which is precisely
what the Healy-Driscoll algorithm is designed to do. To introduce this procedure, we first note that, from
(\ref{57}), 
\begin{equation}\label{65}
f\cdot p_{l+K}=Q_{l,K}\cdot (f\cdot p_l)+R_{l,K}\cdot(f\cdot p_{l-1}),
\end{equation}
and therefore
\begin{equation}\label{66}
\mathcal{T}_{K}(f\cdot p_{l+K})=\mathcal{T}_{K}(Q_{l,K}\cdot (f\cdot p_l)+R_{l,K}\cdot(f\cdot p_{l-1})).
\end{equation}
On the other hand, from (\ref{58}), 
$$Q_{l,K}T_{2K+i}\in \mbox{span}_{\mathbb R}\{T_{K+i},\dots,T_{3K+i}\},$$ 
$$R_{l,K}T_{2K+i}\in \mbox{span}_{\mathbb R}\{T_{K+i+1},\dots,T_{3K+i-1}\}.$$ 
Therefore,
\begin{multline}\label{67}
\mathcal{T}_{K}(Q_{l,K}\cdot (f\cdot p_l)+R_{l,K}\cdot(f\cdot p_{l-1}))
= \mathcal{T}_{K}(Q_{l,K}\cdot \mathcal{T}_{2K}(f\cdot p_l)+R_{l,K}
\cdot\mathcal{T}_{2K}(f\cdot p_{l-1})). 
\end{multline}

The degree of $Q_{l,K}\cdot \mathcal{T}_{2K}(f\cdot p_l)$ is less than $3K$ and by applying the Lagrange 
interpolation operator, we obtain
$$
 Q_{l,K}\cdot \mathcal{T}_{2K}(f\cdot p_l)-\mathcal{S}_{2K}(Q_{l,K}\cdot \mathcal{T}_{2K}(f\cdot p_l)) 
  =T_{2K}\cdot w(x),
$$
where $w(x)$ is polynomial of degree less than $K$ and we have used the fact that the Chebyshev points 
$x_j^{2K}$ are the zeros of $T_{2K}$.
 
Again, from (\ref{58}), 
$$
T_{2K}\cdot w(x)\in \mbox{span}_{\mathbb R}\{T_{K},\dots,T_{3K}\}.  
$$
Therefore,
\begin{equation}\label{68}
\mathcal{T}_{K}({T_{2K}\cdot w(x)})=0. 
\end{equation}
Finally, from (\ref{66})-(\ref{68}),
\begin{equation}\label{69}
\begin{split}
Z_{l+K}^K
       & =\mathcal{T}_{K}(Q_{l,K}\cdot (f\cdot p_l)+R_{l,K}\cdot(f\cdot p_{l-1}))\\
       & =\mathcal{T}_{K}(Q_{l,K}\cdot \mathcal{T}_{2K}(f\cdot p_l)+R_{l,K}
           \cdot\mathcal{T}_{2K}(f\cdot p_{l-1}))  \\
       & =\mathcal{T}_{K}(\mathcal{S}_{2K}(Q_{l,K}\cdot \mathcal{T}_{2K}(f\cdot p_l)+
           R_{l,K}\cdot\mathcal{T}_{2K}(f\cdot p_{l-1})))\\  
       & =\mathcal{T}_{K}(\mathcal{S}_{2K}(Z_{l}^{2K}\cdot Q_{l,K})
          +\mathcal{S}_{2K}(Z_{l-1}^{2K}\cdot R_{l,K})).\\  
\end{split}
\end{equation}

With a similar argument, we can also show that
\begin{equation}\label{70}
  Z_{l+K-1}^K= \mathcal{T}_{K}(\mathcal{S}_{2K}(Z_{l}^{2K}\cdot Q_{l,K-1})
+\mathcal{S}_{2K}(Z_{l-1}^{2K}\cdot R_{l,K-1})).
\end{equation}

The Healy-Driscoll algorithm is based on formulas (\ref{69}) and (\ref{70}). From (\ref{63}),
$Z_{l+K-1}^K$, $Z_{l+K}^K$ can be computed from $Z_{l-1}^{2K}$ and $Z_{l}^{2K}$ using the FCT in
$O(K\log K)$ operations.
At stage 0, the algorithm computes ($Z_0^j,\mbox{  }Z_1^j$) using FCT.
At stage 1, we get ($Z_{N/2}^{N/2},\mbox{  }Z_{N/2+1}^{N/2}$) from ($Z_{0}^{N},\mbox{  }Z_{1}^{N}$) 
using (\ref{69}),(\ref{70}). At stage 2, we compute ($Z_{N/4}^{N/4},\mbox{  }Z_{N/4+1}^{N/4}$) from 
($Z_{0}^{N/2},\mbox{  }Z_{1}^{N/2}$) and ($Z_{3/4N}^{N/4},\mbox{  }Z_{3/4N+1}^{N/4}$) from 
($Z_{N/2}^{N/2},\mbox{  }Z_{N/2+1}^{N/2}$), etc; see figure 2.
The total number of stages is $\log_2(N)$ and each stage costs $O(N\log N)$ operations, 
for an overall cost of $O(N\log^2 N)$ operations. 

\epsfxsize=10 cm
\epsfysize=10 cm
\begin{center}
\epsfbox{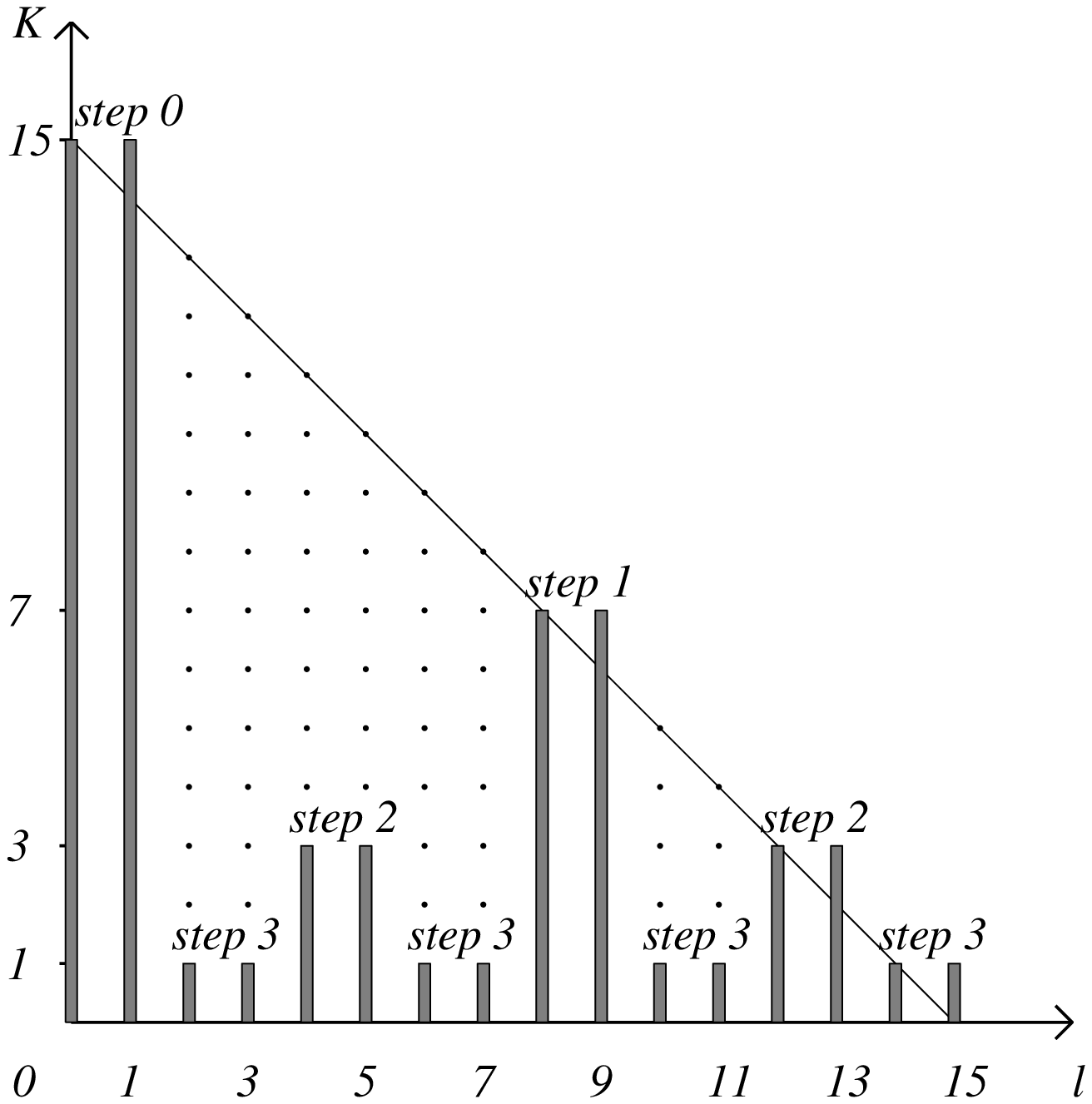}
\end{center}
$$
  \mbox{Figure 2: Computation of } Z_l^K \mbox{ with } N=16
$$

\subsubsection{The Inverse Fast Legendre Transform (IFLT)}
Clearly the IFLT can be written as a matrix in the form
\begin{equation*}
\mathcal{M}_{N}=
\left(
\begin{matrix}
p_{0}(x_{0}^{N}) & p_{1}(x_{0}^{N})&\dots & p_{N-1}(x_{0}^{N})\\
p_{0}(x_{1}^{N}) & p_{1}(x_{1}^{N})&\dots & p_{N-1}(x_{1}^{N})\\
\dots & \dots & \dots &\dots\\
p_{0}(x_{N-1}^{N}) & p_{1}(x_{N-1}^{N})&\dots & p_{N-1}(x_{N-1}^{N})\\
\end{matrix}
\right).
\end{equation*}
As shown in \cite{healy2},
\begin{equation*}
\mathcal{M}_{N}^t=
\left(
\begin{matrix}
p_{0}(x_{0}^{N}) & p_{0}(x_{1}^{N})&\dots & p_{0}(x_{N-1}^{N})\\
p_{1}(x_{0}^{N}) & p_{1}(x_{1}^{N})&\dots & p_{1}(x_{N-1}^{N})\\
\dots & \dots & \dots &\dots\\
p_{N-1}(x_{0}^{N}) & p_{N-1}(x_{1}^{N})&\dots & p_{N-1}(x_{N-1}^{N})\\
\end{matrix}
\right)
\end{equation*}
can be applied to a vector in $O(N(\log N )^2)$ operations. Indeed  
$\mathcal{M}_{N}^t$ can be decomposed into $\log_2 N$ factors, 
\begin{equation}\label{71}
\mathcal{M}_{N}^t=A_{\log_2 N} A_{\log_2 N-1} \cdot \cdot \cdot A_{2} A_{1}
\end{equation}
and $A_{i}^t$ can itself be applied in $O(N(\log N ))$ operations.
Therefore we have that
\begin{equation}\label{72}
IFLT_{N}=\mathcal{M}_{N}=
A_{1}^{t} A_{2}^{t} \cdot \cdot \cdot A_{\log_2 N-1}^{t} A_{\log_2 N}^{t}
\end{equation} 
can be applied in $O(N\log^2 N)$ operations; for details, see \cite{healy2}.

\begin{lemma}  {IFLT algorithm.}\\
\textbf{Input} \mbox{  }
$\mathbf{\hat{f}}=(\hat{f}_{0},\dots,\hat{f}_{N-1})$: Legendre transform coefficients and $N$ is a power of 2.\\
\textbf{Output} 
$\mathbf{f}=(f_{0},\dots,f_{N-1})$: Inverse Legendre transform values.\\
\textbf{Stages}\\
0 stage,
$$\mathbf{a}=[\hat{f}_{0},0,\hat{f}_{1},0,\dots,\hat{f}_{N-1},0]$$
where $a_{i}=[\hat{f}_{i},0] \mbox{  and  }\mathbf{a}=[a_{0},\dots,a_{N-1}]$\\
k stage,$\mbox{ for k=1 to}\log_{2}{N}-1 $
\renewcommand{\labelenumi}{{\rm(\alph{enumi})}}
\begin{enumerate}
$K=2^k$\\
\item calculate $\mathbf{c}\in \mathcal{M}_{1*4N}$,\\
$$\mathbf{c}=[c_{0},\dots,c_{2N/K-1}],\mbox{  }c_{i}=\mbox{ cheby }^{-1}_{2K}([a_{i},o_{1*K}]),\mbox{  }c_{i}\in \mathcal{M}_{1*2K}$$
where
$o_{1*K}:\mbox{ zero matrix of size } 1*K$\\
\item calculate $\mathbf{g}\in \mathcal{M}_{1*2N}$,
\begin{equation*}
\left(
\begin{matrix}
A^K_{1} & O & \dots & O\\
O & A^K_{2}& \dots & O\\
\\
O & O & \dots & A^K_{N/(2K)}
\end{matrix}
\right)
\left(
\begin{matrix}
c_{0}^{t} \\
c_{1}^{t} \\
. \\
. \\
. \\
c_{2N/K-1}^{t}
\end{matrix}
\right)=
\left(
\begin{matrix}
g_{0}^t \\
g_{1}^t \\
. \\
g_{N/K-1}^t
\end{matrix}
\right),
\mbox{  }g_{i}\in \mathcal{M}_{1*2K}
\end{equation*}
where
\begin{equation*}
A_i^{K}=
\left(
\begin{matrix}
I_{2K*2K} & o_{2K*2K} & R_{l}^{K-1} & R_{l}^{K}\\
o_{2K*2K} & I_{2K*2K} & Q_{l}^{K-1} & Q_{l}^{K}
\end{matrix}
\right)\in \mathcal{M}_{4K*8K}
\end{equation*}

$$ 
R_{l}^{K-1}=\mbox{ diag }(R_{l}^{K-1}(x_{0}^{2K}),\dots,R_{l}^{K-1}(x_{2K-1}^{2K}))\in \mathcal{M}_{2K*2K}
$$
$$
R_{l}^{K}=\mbox{ diag }(R_{l}^{K}(x_{0}^{2K}),\dots,R_{l}^{K}(x_{2K-1}^{2K}))
$$
$$
Q_{l}^{K-1}=\mbox{ diag }(Q_{l}^{K-1}(x_{0}^{2K}),\dots,Q_{l}^{K-1}(x_{2K-1}^{2K}))
$$
$$
Q_{l}^{K}=\mbox{ diag }(Q_{l}^{K}(x_{0}^{2K}),\dots,Q_{l}^{K}(x_{2K-1}^{2K}))
$$
$$
  l=2K*(i-1)+1
$$

\item calculate $\mathbf{a}\in \mathcal{M}_{1*2N}$,\\

$$\mathbf{a}=[a_{0},\dots,a_{N/K-1}],\mbox{  }a_{i}=\mbox{ cheby }_{2K}(g_{i})$$
\end{enumerate}
$\log_2{N}$ stage,\\

\begin{itemize}
\item Given $\mathbf{a}=[a_{0},a_{1}]$, compute $\mathbf{e}\in \mathcal{M}_{1*2N}$\\

$$\mathbf{e}=[\mbox{ icheby}_{N}(a_{0}), \mbox{ icheby}_{N}(a_{1})]$$

\item $\mathbf{f}=(f_{0},\dots,f_{N-1})$ is computed as

\begin{equation*}
\left(
\begin{matrix}
diag(P_0(x_{j}^N)) & diag(P_1(x_{j}^N))
\end{matrix}
\right)
\left(
\begin{matrix}
\mathbf{e}^t
\end{matrix}
\right)=
\left(
\begin{matrix}
\mathbf{f}^{t}
\end{matrix}
\right)
\end{equation*}
\end{itemize}
\end{lemma}
\subsubsection{Numerical examples }
In table 5 and figure 3 we present results representative of the performance of the algorithms described in
A.1.3 and A.1.4.
To assess the stability and accuracy of the algorithms we shall compare their outcome to \emph{direct} 
calculations of the transform and its inverse. To this end, we let
\begin{equation}\label{73}
DLT(\{ f(x^{2N}_i)\}_{i=0}^{2N-1})=\{c_m\}_{m=0}^{N-1}
\end{equation}
denote the forward transform, and 
\begin{equation}\label{74}
IDLT(\{c_m\}_{m=0}^{N-1})=\{ f(x^{2N}_i)\}_{i=0}^{2N-1}
\end{equation}
the inverse, computed directly, without acceleration and using the three-term recurrence (\ref{55})
to evaluate Legendre polynomials. We further introduce the notation $DLT^{QP}$ and $IDLT^{QP}$ to denote 
the same transforms but with the three-term recurrence (\ref{55}) evaluated in quadruple precision. 
In table 5 we compare the errors
$$
\|DLT_{N}(IDLT_{N}^{QP}(v_N^{N/2}))-v_N^{N/2}\|_{\infty},
$$ 
$$
\|IDLT_{N}(v_N^{N/2})-IDLT_{N}^{QP}(v_N^{N/2})\|_{\infty}
$$ 
labeled ``DLT'' and ``IDLT'' respectively, with that incurred by the fast transforms, namely
$$
\|FLT(IDLT_{N}^{QP}(v_N^{N/2}))-v_{N}^{N/2}\|_{\infty}\ \mbox{    (``FLT'')}
$$ 
$$
\|IFLT(v_{N}^{N/2})-IDLT_{N}^{QP}(v_N^{N/2})\|_{\infty}\ \mbox{    (``IFLT'')}  
$$
where $v_{m}^{n}$ is 
a vector of size $m$ with unit $n$th component and zeros everywhere else. 
In addition, the table includes results obtained on application of the fast transforms when the 
recurrence (\ref{56}) are precomputed in quadruple precision (``FLT-QP'' and ``IFLT-QP''). 
We see that this latter 
strategy provides errors that are comparable to double precision DLT and IDLT calculations.
Finally, in figure 3, we display timings for the FLT algorithm showing that, indeed, the operation count 
is proportional to $N\log^2 N$.\\

\begin{center}
\renewcommand{\arraystretch}{1.25}
\begin{tabular}{|r| p{.7in}   | p{.7in}   | p{.7in}  |  p{.7in} | p{.7in} | p{.7in} |}
\hline
 \textbf{N}  & \textbf{DLT} & \textbf{FLT} & \textbf{FLT-QP}& \textbf{IDLT}  &\textbf{IFLT} &\textbf{IFLT-QP} \\ \hline
256    & 1.98476e-14  &2.8107e-13  & 1.17584e-14& 2.52909e-13  & 5.27689e-13   & 1.99396e-13    \\ \hline 
1024   & 8.0331e-14  &3.49153e-12 & 3.21525e-14& 4.18487e-12  & 2.3919e-10  & 4.22617e-12  \\ \hline
4096  & 1.49405e-13 &1.02863e-10 & 1.51153e-13 & 8.09853e-11& 1.45763e-08  & 3.95126e-11  \\ \hline
16384   & 9.05882e-13  &1.27641e-09 & 5.2725e-13  & 1.2968e-09& 4.37189e-07   & 5.39e-10   \\ \hline
\end{tabular}
\end{center}
Table 5: Examples of errors in the FLT and IFLT algorithms, QP indicates that 
the (pre-) computation  of the recurrences (\ref{56}) are carried out in 
quadruple precision. \\

\epsfxsize=7 cm
\epsfysize=5 cm
\begin{center}
\epsfbox{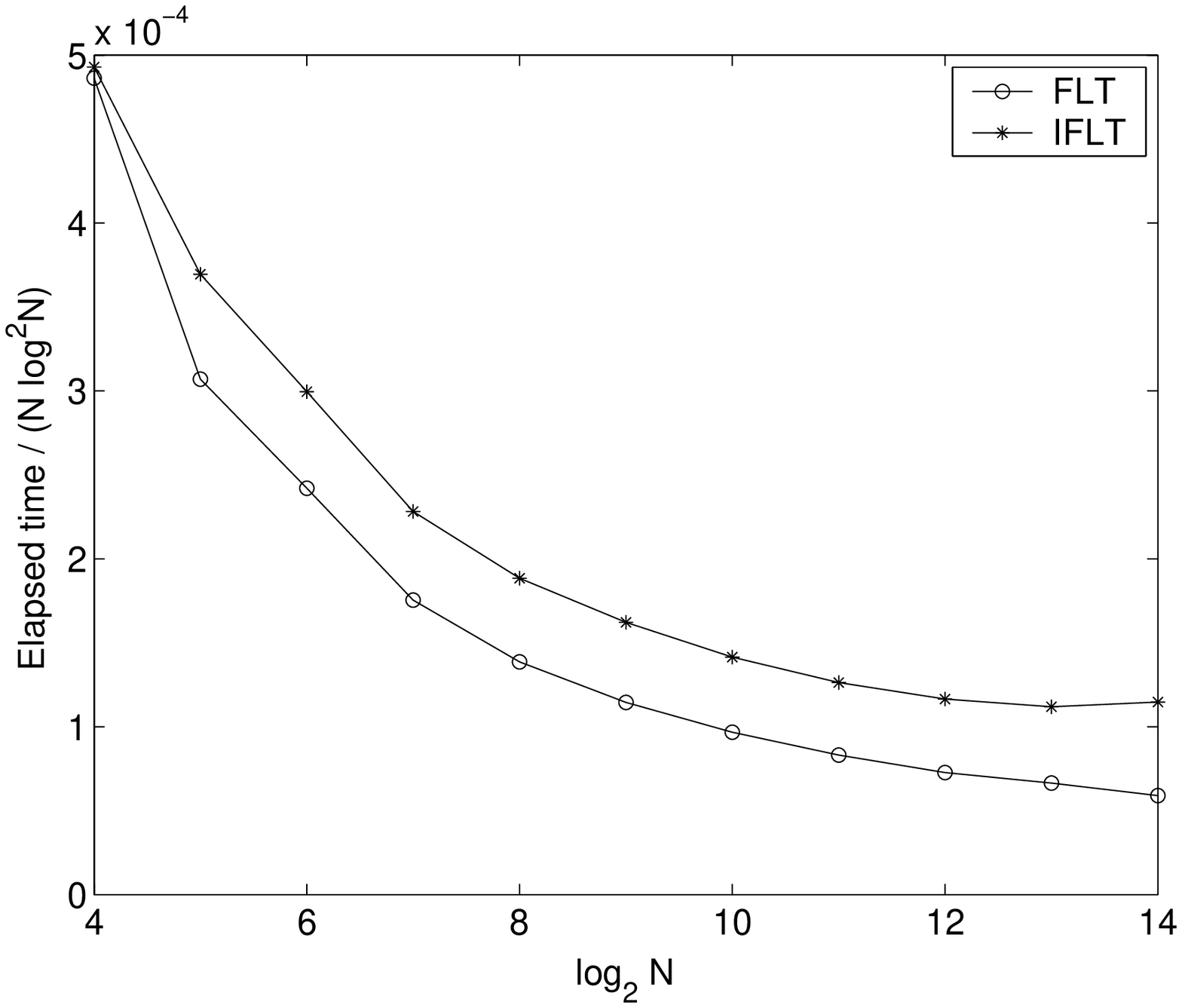}
\end{center}
$$
  \mbox{Figure 3: Timings for the evaluation of FLT}(\{ f(x^{2N}_i)\}_{i=0}^{2N-1}) 
$$
$$
\mbox{ and IFLT}(v_{N}^{N/2}) \mbox{ for different values of N.}
$$

\subsection{Computation of modified Bessel functions}
The spherical Bessel functions $j_{\nu}(\rho),  y_{\nu}(\rho), h_{\nu}^{(1)}(\rho), h_{\nu}^{(2)}(\rho)$ are defined in terms of ordinary Bessel functions by the relations,
\begin{equation}\label{75}
j_{\nu}(\rho)=\sqrt{\frac{\pi}{2\rho}}J_{\nu+\frac{1}{2}}(\rho),
\end{equation}
\begin{equation}\label{76}
y_{\nu}(\rho)=\sqrt{\frac{\pi}{2\rho}}Y_{\nu+\frac{1}{2}}(\rho),
\end{equation}
\begin{equation}\label{77}
h_{\nu}^{(1),(2)}(\rho)=\sqrt{\frac{\pi}{2\rho}}H_{\nu+\frac{1}{2}}^{(1),(2)}(\rho),
\end{equation}
where
\begin{equation}\label{78}
h_{\nu}^{(1),(2)}(\rho)= j_{\nu}(\rho)\pm iy_{\nu}(\rho).
\end{equation}
A standard calculation of Bessel functions evaluates these from the recurrence
\begin{equation}\label{83}
 Z_{\nu+1}(\rho)=
\frac{2\nu}{\rho}Z_{\nu}(\rho)- Z_{\nu-1}(\rho), 
\end{equation}
$$
 Z_{\nu}\in\{j_{\nu},\mbox{  } y_{\nu},\mbox{  } h^{(1),(2)}_{\nu}\}.
$$
Our modified Bessel functions are as follows.
\begin{equation}\label{84}
\widetilde{j_n}(\rho)=\frac{1\cdot3\cdot5\cdot\dots(2n+1)}{\rho^n}j_n(\rho)=
\bigg[1-\frac{\frac{1}{2}\rho^2}{1!(2n+3)}+\frac{({\frac{1}{2}\rho^2})^2}{2!(2n+3)(2n+5)}+\dots
\bigg]
\end{equation}
\begin{equation}\label{85}
\widetilde{y_n}(\rho)=\frac{\rho^{n+1}}{-1\cdot1\cdot3\cdot5\cdot\dots(2n-1)}y_n(\rho)=
\bigg[1-\frac{\frac{1}{2}\rho^2}{1!(1-2n)}+\frac{({\frac{1}{2}\rho^2})^2}{2!(1-2n)(3-2n)}+\dots
\bigg]
\end{equation}
Therefore the new recurrences of these are
\begin{equation}\label{86}
\widetilde{j}_{n-1}(\rho)=\widetilde{j}_n(\rho)-\frac{\rho^2}{(2n+1)(2n+3)}\widetilde{j}_{n+1}(\rho)
\end{equation}
\begin{equation}\label{87}
\widetilde{y}_{n+1}(\rho)=\widetilde{y}_n(\rho)-\frac{\rho^2}{(2n-1)(2n+1)}\widetilde{y}_{n-1}(\rho)
\end{equation}
where,
\begin{equation*}
\begin{split}
& \widetilde{y}_0(\rho)=-\rho\cdot\sqrt{\frac{\pi}{2\rho}}Y_{\frac{1}{2}}(\rho)\\ 
& \widetilde{y}_1(\rho)=-\rho^2\cdot\sqrt{\frac{\pi}{2\rho}}Y_{\frac{3}{2}}(\rho)\\ 
\end{split}
\end{equation*}
The reason that we need to compute $\widetilde{j}_n(\rho)$ with
downward recurrence is that (\ref{86}) has 2 linearly independent
solutions $\widetilde{j}_n(\rho)$ and
$\breve{y}_n(\rho)=\frac{1\cdot3\cdot5\cdot\dots(2n+1)}{\rho^n}y_n(\rho)$.
Since $\breve{y}_n(\rho)$ is an exponentially growing solution as $n$ increases, the
upward recurrence relation is numerically unstable since the round-error
 in the $\widetilde{j}_0(\rho)$ and  $\widetilde{j}_1(\rho)$ is rapidly
 amplified by the recurrence. On the other hand, the downward recurrence
 is stable since the round-off error in the starting values is quickly
 damped by the recurrence. Computing  $\widetilde{j}_n(\rho)$ and
 $\widetilde{j}_{n-1}(\rho)$ can be done directly from the (\ref{84}).

\end{document}